\documentclass[11pt,a4paper]{amsart}
\usepackage{tikz}
\usepackage[margin=2.6cm]{geometry}
\usepackage[utf8]{inputenc}
\usepackage{amsmath,amssymb,amsthm}
\usepackage[shortlabels]{enumitem}
\usepackage{xcolor}
\usepackage[colorlinks,pagebackref,linkcolor=blue!60!black,citecolor=blue!60!black,urlcolor=blue!60!black]{hyperref}
\renewcommand*{\backref}[1]{}
\renewcommand*{\backrefalt}[4]{\ifcase #1 \or {\small(Cited on page~#2.)}\else {\small(Cited on pages~#2.)}\fi}
\makeatletter
\newcounter{BRcite}
\let\BR@orig@Hy@backout\Hy@backout
\def\Hy@backout#1{\stepcounter{BRcite}\ifhmode\vadjust pre{\hypertarget{brcite.\the\value{BRcite}}{}}\else\hypertarget{brcite.\the\value{BRcite}}{}\fi%
  \begingroup\edef\@currentHref{brcite.\the\value{BRcite}}\def\@currentlabel{cite}%
  \BR@orig@Hy@backout{#1}\endgroup}

\makeatother

\newcommand{\Q}{\mathbb{Q}}
\newcommand{\R}{\mathbb{R}}
\newcommand{\C}{\mathbb{C}}
\newcommand{\Z}{\mathbb{Z}}
\newcommand{\OK}{\mathcal{O}_K}
\newcommand{\Kt}{\widetilde K}
\newcommand{\iu}{\mathrm{i}}
\newcommand{\cD}{\mathcal{D}}
\newcommand{\fB}{\mathfrak{B}}
\DeclareMathOperator{\Gal}{Gal}
\DeclareMathOperator{\Hom}{Hom}
\DeclareMathOperator{\Alt}{Alt}
\DeclareMathOperator{\im}{Im}
\DeclareMathOperator{\re}{Re}
\DeclareMathOperator{\Tr}{Tr}
\newcommand{\NS}{\mathop{\textup{N\'eron--Severi}}\nolimits}

\theoremstyle{plain}
\newtheorem{theorem}{Theorem}[section]
\newtheorem{proposition}[theorem]{Proposition}
\newtheorem{lemma}[theorem]{Lemma}
\newtheorem{corollary}[theorem]{Corollary}
\theoremstyle{definition}
\newtheorem{definition}[theorem]{Definition}
\theoremstyle{remark}
\newtheorem{remark}[theorem]{Remark}
\newtheorem{example}[theorem]{Example}

\begin{document}

\title[Quasi-abelian Cousin groups of number fields]{Quasi-abelian Cousin groups of number fields:\\ a Galois-orbit criterion and the complete list}
\author{Karl Oeljeklaus}
\address{Aix Marseille Univ, CNRS, I2M, 3 place Victor Hugo, Case 19, 13331 Marseille Cedex 3, France}
\email{karl.oeljeklaus@univ-amu.fr}
\date{}
\dedicatory{Dedicated to the memory of Aldo Andreotti and Alan T.~Huckleberry}

\begin{abstract}
Let $K$ be a number field with at least one complex place. A choice $S$ of one embedding $K\to\C$ at each place
embeds the ring of integers of $K$ as a discrete subgroup of $\C^{n}$, where $n$ is the number of places
of $K$. The quotient $T_{K,S}$ is a
Cousin group, or a complex torus if $K$ is totally imaginary. We determine exactly when $T_{K,S}$ is quasi-abelian. The answer is a Galois criterion: no
complex conjugation of the Galois closure of $K$ may interchange two of the chosen complex embeddings.
Equivalently, $S$ is lifted from a subfield with at most one complex place above each place of the largest
totally real subfield of $K$. For totally imaginary $K$ this is the theorem of Shimura and Taniyama on
abelian varieties with complex multiplication; our results extend it to all number fields. We also obtain
the rank of the N\'eron--Severi group, and show that $T_{K,S}$ is either quasi-abelian or has only constant
meromorphic functions.

\medskip
\noindent{\scshape R\'esum\'e.} \emph{Groupes de Cousin quasi-ab\'eliens des corps de nombres~: un crit\`ere
d'orbites galoisiennes et la liste compl\`ete.}
Soit $K$ un corps de nombres ayant au moins une place complexe. Le choix $S$ d'un plongement $K\to\C$ en chaque place
plonge l'anneau des entiers de $K$ comme sous-groupe discret de $\C^{n}$, o\`u $n$ est le nombre de places
de $K$. Le quotient $T_{K,S}$ est un
groupe de Cousin, ou un tore complexe si $K$ est totalement imaginaire. Nous d\'eterminons exactement quand $T_{K,S}$ est quasi-ab\'elien. La r\'eponse est un
crit\`ere galoisien~: aucune conjugaison complexe de la cl\^oture galoisienne de $K$ ne doit \'echanger deux
des plongements complexes choisis. De mani\`ere \'equivalente, $S$ est relev\'e d'un sous-corps ayant au plus
une place complexe au-dessus de chaque place du plus grand sous-corps totalement r\'eel de $K$. Pour $K$
totalement imaginaire, c'est le th\'eor\`eme de Shimura et Taniyama sur les vari\'et\'es ab\'eliennes \`a
multiplication complexe~; nos r\'esultats l'\'etendent \`a tous les corps de nombres. Nous obtenons aussi le
rang du groupe de N\'eron--Severi, et montrons que $T_{K,S}$ est soit quasi-ab\'elien, soit d\'epourvu de
fonctions m\'eromorphes non constantes.
\end{abstract}

\maketitle

\section{Introduction}\label{sec:intro}

Let $K$ be a number field with $r_1$ real places and $r_2$ complex places, so that its degree is
$\delta=r_1+2r_2$, and let $\OK$ be its ring of integers. Choose one embedding $K\to\C$ at each place; we call
this choice $S$. Map $\OK$ into $\C^{r_1+r_2}$ with these embeddings, and divide out. The result, $T_{K,S}$, is
a complex abelian Lie group of dimension $r_1+r_2$. When $K$ is totally imaginary ($r_1=0$) it is a complex
torus, and the question of
when it is an abelian variety was settled by Shimura and Taniyama~\cite{ST}. When $K$ has real and complex
places, $T_{K,S}$ is a Cousin group: it is not compact and has no nonconstant holomorphic
functions~\cite{AG73,OT,AbeFM}; see Theorem~\ref{thm:mero}. The right analogue of an abelian variety is then a
\emph{quasi-abelian} group in the sense of Andreotti and Gherardelli, one that carries an ample Riemann form.
The name goes back to Severi. For Cousin groups, being quasi-abelian is the same as being
quasi-projective; this was proved in general by Abe~\cite{Abe89}, see~\cite[Main Thm.~3.2.21]{AK}.

The study of these groups goes back to Andreotti and Gherardelli~\cite{AG73,AG}; see
also~\cite{Gher,Dod80,AbeFM,Um17}. They are natural objects. Dodson showed that a quasi-abelian
variety whose endomorphism algebra contains a number field of the largest possible degree comes from the
embeddings of that field, with a lattice of the field in place of $\OK$~\cite{Dod80}. The same groups occur in
complex geometry. Let $\mathbb H$ be the upper half-plane. Then $T_{K,S}$ contains the domain
$(\mathbb H^{r_1}\times\C^{r_2})/\OK$, and the compact quotients of this domain by suitable groups of units are
the manifolds of Oeljeklaus and Toma~\cite{OT,BO}.

Another motivation comes from the classical correspondence between number fields and commutative algebraic
groups. The field $\Q$ gives $\C/\Z\cong\C^{\times}$, whose torsion points, the roots of unity, generate the
abelian extensions of $\Q$. Imaginary quadratic fields give elliptic curves with complex multiplication, and CM fields give the
abelian varieties of Shimura and Taniyama. Valli\`eres proposed to continue this list with the groups
$T_{K,S}$ of an arbitrary number field, with $\OK$ replaced by any lattice of $K$. He
studied in detail the cubic fields that are not totally real; they give Cousin groups of dimension $2$ and
rank $3$~\cite{Val}. For $r_2=1$, Dioguardi Burgio, Falcone and Galici showed that $T_{K,S}$ is an extension of
an elliptic curve $C$ by $(\C^{\times})^{r_1}$. They represented it as a generalized Jacobian of $C$, and
described its torsion points through fractional ideals of $K$ in the cubic case~\cite{DFG}. The periods of the
generalized Jacobians of an elliptic curve had been determined by Di Bartolo and Falcone~\cite{DBF}; see
also~\cite{DF,DB}. Under assumptions on the discriminant, Dioguardi Burgio showed that every number field with
$r_2=1$ arises as the endomorphism algebra of a toroidal group whose real rank exceeds its complex rank by
one~\cite{DB2}. Generalized Jacobians go back to Rosenlicht~\cite{Ros}. They are commutative algebraic
groups that extend the Jacobian of a curve by a linear group, and Serre made them the basis of his geometric
class field theory~\cite{Se}. For $r_2\ge2$ the natural analogue is an extension of an abelian variety by a
torus. Such extensions are exactly the quasi-abelian groups of kind $0$ in the sense of Andreotti and
Gherardelli~\cite{AG73,AK}. This is one
reason to ask which groups $T_{K,S}$ are quasi-abelian at all.

Which pairs $(K,S)$ give quasi-abelian groups? For $r_2=1$ all of them do~\cite{AG73,Gher}. Several families
are known~\cite{AG73,AG,Dod79,Abe13}, and they all come from two constructions. The first, which we call the \emph{trace-form criterion}, uses a
totally real subfield $K_0$ of $K$ above whose places the complex places of $K$ are spread out: at most one above each.
The name recalls how these groups were first shown to be quasi-abelian: by a Riemann form of the shape
$\Tr_{K_0/\Q}\circ\Phi$, the trace from $K_0$ to $\Q$ of a suitable $K_0$-valued alternating form $\Phi$ on $K$.
For a CM field and $K_0=K^{+}$ this is the classical form $\Tr_{K/\Q}(\nu x\bar y)$ of Shimura and Taniyama,
with $\bar\nu=-\nu$, since $\Tr_{K/K_0}(\nu x\bar y)=\nu(x\bar y-\bar xy)$ lies in $K_0$.
The second, due to Dodson, \emph{lifts} a quasi-abelian group from a subfield~\cite{Dod79}. Dodson also
proved a first converse. If $T_{K,S}$ has an ample Riemann form of a special, diagonal shape, then $S$ is
lifted from a subfield that is quadratic over a totally real field~\cite[Prop.~2.2]{Dod79}. Abe showed that
every $T_{K,S}$ carries a Hermitian form that is positive in the complex directions and whose imaginary part
takes values in the ring of integers of a real number field, not necessarily in
$\Z$~\cite[Thm.~4.3]{AbeFM}. He remarked that it was not known how far such groups are from being
quasi-abelian. In the other direction, Umeno proved that for $K=\Q(2^{1/5})$ and one choice of $S$ the N\'eron--Severi
group of $T_{K,S}$ vanishes, so that $T_{K,S}$ has only constant meromorphic functions~\cite{Um17}. Dioguardi
Burgio extended this, by explicit computation with Umeno's period relations, to the pure monogenic quintic
fields and to a class of quintic fields with an essential polynomial, for one choice of $S$ in each
case~\cite[\S3.2.1]{DB}. In this
paper we answer the question completely. The answer comes in two forms.

The first is a criterion in terms of the Galois group (Theorem~\ref{thm:main}): $T_{K,S}$ is quasi-abelian if
and only if no complex conjugation of the Galois closure of $K$ interchanges two of the chosen complex
embeddings. Here a complex conjugation is the automorphism induced by complex conjugation under some
embedding of the Galois closure into $\C$. The idea behind it is simple. Write a rational alternating form on $K$ in coordinates indexed by
the embeddings. Being the imaginary part of a Hermitian form makes some coefficients vanish, while positivity
requires others to be nonzero. But the Galois group permutes the coefficients of a \emph{rational} form, so a
zero spreads along a whole orbit. That gives one direction. For the other we show that rational forms are
dense among the real forms with the required zeros, and take one close to a real form that is positive by
construction.

The second form is a list (Theorem~\ref{thm:list}). The two constructions above, used together, give
everything: $T_{K,S}$ is quasi-abelian if and only if $S$ is lifted from a subfield $K'$ such that above each
place of the largest totally real subfield $K^{+}$ of $K$ there is at most one complex place of $K'$. This
subfield is unique, and it does not depend on $S$. By Dodson's theorem, every quasi-abelian $T_{K,S}$ is then
isogenous to a power of a group to which the trace-form criterion applies directly
(Corollary~\ref{cor:isogeny}). For a given $K$ the number of quasi-abelian groups $T_{K,S}$ is
therefore $0$ or $2^{r_2(K_c)}$. Here $K_c$ is a subfield determined by the Galois group and complex
conjugation alone, and $r_2(K_c)$ is its number of complex places. The proof starts from the first criterion and does not come back to Riemann forms.
Everything happens in one graph $\Gamma$. Its vertices are the embeddings of $K$, and two of them are joined
when some complex conjugation interchanges them. A quasi-abelian $S$ gives a proper $2$-colouring of the complex
embeddings: two embeddings of the same colour are never joined. This turns out to be a strong restriction: each connected component of $\Gamma$ is then
bipartite or complete multipartite, and the parts form a block system for the Galois group. The subfield
$K_c$ is the one that belongs to this block system. The few notions from graph theory that we use are recalled in
Sections~\ref{sec:graph} and~\ref{sec:canonical}.

The remaining sections draw consequences. The rank of the N\'eron--Severi group is a sum of orbit lengths
(Section~\ref{sec:NS}). The trace-form criterion and Dodson's lifting theorem get short proofs from the orbit
criterion (Section~\ref{sec:subfields}). The theorem of Shimura and Taniyama follows in two ways: as the case
$r_1=0$ of the list, and by a direct counting argument (Section~\ref{sec:consequences}). There is also a
converse of the trace-form criterion. It says that $T_{K,S}$ is quasi-abelian for \emph{every} $S$ exactly
when the trace-form criterion applies with $K_0=K^{+}$. A field without subfields other than $\Q$ and itself,
for instance a field of prime degree, gives quasi-abelian groups only when $r_2=1$. If the criterion is not
satisfied, then $T_{K,S}$ has only constant meromorphic functions (Section~\ref{sec:meromorphic}). Section~\ref{sec:examples} contains examples.

\section{Setting and statements}\label{sec:setting}

\subsection{The groups \texorpdfstring{$T_{K,S}$}{T(K,S)}}
Let $\Sigma=\Hom(K,\C)$. It consists of the real embeddings $\sigma_1,\dots,\sigma_{r_1}$ and of $r_2$
pairs $\{\tau_j,\bar\tau_j\}$ of complex conjugate ones. We always assume $r_2\ge1$. A \emph{system of
representatives} is a set
\[
  S=\{\sigma_1,\dots,\sigma_{r_1},\tau_1,\dots,\tau_{r_2}\}
\]
containing one embedding from each place: all the real embeddings, and one embedding from each pair of
complex conjugate ones. There are $2^{r_2}$ such systems. Once $S$ is fixed, we label the complex
embeddings so that $\tau_1,\dots,\tau_{r_2}$ are those in $S$; their conjugates $\bar\tau_1,\dots,\bar\tau_{r_2}$
are the complex embeddings not in $S$. Put $g=r_1+r_2$ and
\[
  u\colon K\to\C^{g},\quad u(a)=\bigl(\rho(a)\bigr)_{\rho\in S},\qquad \Lambda=u(\OK),\qquad
  T_{K,S}=\C^{g}/\Lambda .
\]
The $\R$-linear extension $u_\R\colon K\otimes_\Q\R\to\R^{r_1}\times\C^{r_2}$ is an isomorphism, so
$\Lambda$ is discrete of rank $\delta$, its real span is $\R\Lambda=\R^{r_1}\times\C^{r_2}$, and the largest
complex subspace of $\R\Lambda$ is
\[
  F=\R\Lambda\cap\iu\,\R\Lambda=\{0\}\times\C^{r_2}.
\]

\subsection{Quasi-abelian groups}
Following Andreotti and Gherardelli~\cite{AG,AK}, an \emph{ample Riemann form} for $T_{K,S}$ is a Hermitian form
$\widetilde H$ on $\C^{g}$ such that
\[
  \widetilde H \text{ is positive definite on } F \qquad\text{and}\qquad
  \im\widetilde H(\Lambda\times\Lambda)\subseteq\Z ,
\]
and $T_{K,S}$ is \emph{quasi-abelian} if it has one. Positivity is asked for only on $F$, that is, only in
the $r_2$ coordinates of the complex places. No condition is imposed in the real directions; this is
what makes the non-compact case different. For $r_1=0$ we have $F=\C^{g}$, and an ample Riemann form is a
Riemann form in the classical sense.

Andreotti and Gherardelli made this definition to describe the quotients that carry the structure of an
algebraic variety~\cite{AG73,AG}. That a toroidal group is quasi-abelian if and only if it is quasi-projective
was proved in general by Abe~\cite{Abe89}; see~\cite[Thm.~3.1.18 and Main Thm.~3.2.21]{AK}. For $r_1=0$ this is the classical statement that a complex torus with a Riemann form is
an abelian variety.

\subsection{Galois data}\label{subsec:galois}
Let $\Kt$ be the Galois closure of $K$ in $\C$, and $G=\Gal(\Kt/\Q)$. In general $\Kt$ is strictly larger than $K$. Indeed, $\Kt=K$
means that $K$ is Galois over $\Q$; then all embeddings of $K$ have the same image, so $K$ is totally real or
totally imaginary, and the first case is excluded by $r_2\ge1$. Thus $\Kt\ne K$ whenever $r_1\ge1$. For $r_1=0$,
equality holds exactly for the Galois fields, for instance the cyclotomic fields, but also for non-CM fields such as
$\Q(\sqrt[3]2,\zeta_3)$. We regard $\Kt$ once and for all as a
fixed subfield of $\C$. Every embedding of $K$ into $\C$ then takes its values in $\Kt$, so that
$\Sigma=\Hom(K,\Kt)$. The group $G$ acts on $\Sigma$ by $\gamma\cdot\rho=\gamma\circ\rho$, transitively. Complex
conjugation of $\C$ maps $\Kt$ to itself and restricts to an element $c\in G$. It fixes
the real embeddings and interchanges $\tau_j$ and $\bar\tau_j$. Its conjugates $\gamma c\gamma^{-1}$ are the
complex conjugations that $\Kt$ inherits from its other embeddings into $\C$; we call all of them
\emph{complex conjugations of $\Kt$}. In general they differ from $c$. The letter $c$ is reserved for the
distinguished one, the restriction of complex conjugation of $\C$, and the words \emph{real} and
\emph{complex} for embeddings always refer to it: $\rho$ is real exactly when $c\rho=\rho$. We write $\cD$ for the set of complex conjugations and
$N=\langle\cD\rangle$ for the normal subgroup of $G$ they generate.

The group $G$ also acts on the set of unordered pairs $\{\rho,\rho'\}$ of distinct embeddings. Two kinds of
pairs matter:
\begin{itemize}[leftmargin=0.8cm,itemsep=2pt]
\item the \emph{diagonal pairs} $\{\tau_j,\bar\tau_j\}$, $j=1,\dots,r_2$. They are the $2$-cycles of $c$
  and do not depend on $S$;
\item the \emph{forbidden pairs} $\{\tau_j,\tau_k\}$, $j\ne k$: the pairs of two distinct complex
  embeddings that both belong to $S$. Unlike the diagonal pairs they depend on $S$: if $\bar\tau_j$ is put
  into $S$ in place of $\tau_j$, the forbidden pairs $\{\tau_j,\tau_k\}$ become $\{\bar\tau_j,\tau_k\}$.
  They exist only if $r_2\ge2$. The name comes from Lemma~\ref{lem:dictionary}: a form satisfying (H) must
  vanish on them.
\end{itemize}

\begin{theorem}[orbit criterion]\label{thm:main}
Let $K$ be a number field with $r_2\ge1$, and $S$ a system of representatives. The following are
equivalent.
\begin{enumerate}[(a),leftmargin=0.9cm,itemsep=2pt]
\item $T_{K,S}$ is quasi-abelian.
\item No $G$-orbit of pairs of embeddings contains both a diagonal pair and a forbidden pair.
\item No complex conjugation of $\Kt$ interchanges two distinct complex embeddings $\tau_j$ and $\tau_k$
  that belong to $S$.
\end{enumerate}
\end{theorem}

Conditions (b) and (c) involve only the permutation group $G$ on $\Sigma$, the element $c$ and the set $S$.
Whether $T_{K,S}$ is quasi-abelian is therefore a finite combinatorial question. Three cases can be read off
at once.

\begin{corollary}\label{cor:easy}
\begin{enumerate}[(i),leftmargin=0.9cm,itemsep=2pt]
\item If $r_2=1$, then $T_{K,S}$ is quasi-abelian.
\item If $r_2\ge2$ and $G$ acts transitively on the pairs of embeddings (for instance if $G$ is
  $2$-transitive on $\Sigma$), then $T_{K,S}$ is not quasi-abelian, for any $S$.
\item If $K$ is a CM field, then $c$ is central in $G$, and $T_{K,S}$ is an abelian variety for every $S$.
\end{enumerate}
\end{corollary}

\begin{proof}
For $r_2=1$ there are no forbidden pairs. In case (ii) there is a single orbit, and it contains both kinds
of pairs: diagonal pairs exist because $r_2\ge1$, and forbidden pairs because $r_2\ge2$.

(iii) Write $K=K^{+}(\sqrt{\beta})$ with $K^{+}$ totally real and $\beta\in K^{+}$ totally negative. Then $\Kt$
is generated over the Galois closure $L$ of $K^{+}$, which is totally real, by the square roots of the
conjugates of $\beta$, and these conjugates are totally negative. Under every embedding of $\Kt$ into $\C$,
complex conjugation therefore fixes $L$ and changes the sign of each of these square roots. Hence all complex
conjugations of $\Kt$ coincide, that is, $c$ is central in $G$. The only complex conjugation is then $c$
itself, and it interchanges each $\tau_j$ with $\bar\tau_j\notin S$. So condition (c) of
Theorem~\ref{thm:main} holds, and $T_{K,S}$, which is compact since $r_1=0$, is an abelian variety.
\end{proof}

Part (i) is classical: such a group is fibred over elliptic curves in infinitely many ways~\cite{AG73,Gher};
see also~\cite[Remark~3.2]{Dod79}.

Part (iii) is the classical fact that every CM type of a CM field gives an abelian variety; see
Theorem~\ref{thm:ST} for the general statement in the compact case.

Part (ii) covers fields whose Galois group is the full symmetric group, which is the generic case. It also
covers $K=\Q(2^{1/5})$, of signature $(1,2)$: there $G$ is the Frobenius group of order $20$, which is
$2$-transitive on the five embeddings. This recovers Umeno's result~\cite{Um17}.

\subsection{Subfields, blocks and lifting}\label{subsec:subfields}
The second main result is phrased in terms of subfields. A \emph{subfield} is always a subfield of $K$.

\smallskip\noindent\emph{Embeddings and places.} For a subfield $L$ let $\Sigma(L)=\Hom(L,\C)$. An embedding
$\iota\in\Sigma(L)$ is \emph{real} if $\iota(L)\subseteq\R$, that is, if $c\circ\iota=\iota$. The real places of
$L$ are its real embeddings, and its complex places are the pairs $\{\iota,c\circ\iota\}$ of non-real
embeddings; $r_1(L)$ and $r_2(L)$ denote their numbers. If $L\subseteq L'$ are subfields, a place of $L'$
\emph{lies above} the place of $L$ obtained by restricting its embeddings to $L$.

\smallskip\noindent\emph{Blocks.} Restriction $\rho\mapsto\rho|_L$ maps $\Sigma$ onto $\Sigma(L)$ and commutes
with the action of $G$. Its fibres
\[
  \Sigma_\iota=\{\rho\in\Sigma:\ \rho|_{L}=\iota\},\qquad \iota\in\Sigma(L),
\]
are the \emph{blocks over $L$}. They partition $\Sigma$, each has $[K:L]$ elements, and $G$ permutes them:
$\gamma\Sigma_\iota=\Sigma_{\gamma\iota}$. Conversely, every partition of $\Sigma$ that $G$ permutes arises in
this way from exactly one subfield. (Identify $\Sigma$ with $G/G_K$, where $G_K=\Gal(\Kt/K)$; such partitions
correspond to the subgroups between $G_K$ and $G$, hence to the subfields of $K$.) A larger subfield gives a
finer partition. In particular, the orbits of a normal subgroup of $G$ are the blocks over some subfield.
Since $c\,\Sigma_\iota=\Sigma_{c\iota}$, the embedding $\iota$ is real if and only if $c$ maps the block
$\Sigma_\iota$ to itself; such a block consists of real embeddings and of pairs of complex conjugate ones, and
the latter are the complex places of $K$ above $\iota$.

\smallskip\noindent\emph{The field $K^{+}$.} A subfield is \emph{totally real} if all its embeddings are
real. The compositum of two totally real subfields is again totally real, so $K$ has a largest totally real
subfield, which we denote by $K^{+}$. If a subfield $L$ contains $K^{+}$, then $K^{+}$ is also the largest
totally real subfield of $L$. The blocks over $K^{+}$ have a description in terms of the complex
conjugations.

\begin{lemma}\label{lem:N}
Two embeddings of $K$ lie in the same $N$-orbit if and only if they agree on $K^{+}$; that is, the $N$-orbits
are the blocks over $K^{+}$. A $G$-invariant partition of $\Sigma$ has all its blocks mapped to themselves
by $c$ if and only if it is coarser than the partition into $N$-orbits.
\end{lemma}

\begin{proof}
Let $L$ be a field between $\Q$ and $\Kt$. Its embeddings are the restrictions $\gamma|_L$ with $\gamma\in G$,
and $\gamma|_L$ is real if and only if $\gamma^{-1}c\gamma$ fixes $L$ pointwise. Hence $L$ is totally real if and
only if $N$ fixes $L$. The fixed field of $N$ is therefore the largest totally real subfield of $\Kt$, and its
intersection with $K$ is $K^{+}$. Identify
$\Sigma$ with $G/G_K$. The $N$-orbit of $\gamma G_K$ is $N\gamma G_K=\gamma NG_K$, and $NG_K=\Gal(\Kt/K^{+})$,
because its fixed field is $K\cap\Kt^{N}=K^{+}$. Thus the $N$-orbit of $\gamma|_K$ consists of the embeddings
that agree with $\gamma$ on $K^{+}$. For the last assertion: $c$ maps every block of a $G$-invariant partition to
itself if and only if every conjugate of $c$ does, that is, if and only if $N$ does. This means that every
block is a union of $N$-orbits.
\end{proof}

\smallskip\noindent\emph{Lifting.} Let $K'$ be a subfield and $S'\subseteq\Sigma(K')$ a system of
representatives for $K'$. Following Dodson~\cite{Dod79}, we say that $S$ is \emph{lifted from $S'$} if
\begin{enumerate}[(L1),leftmargin=1.1cm,itemsep=2pt]
\item\label{L1} no complex place of $K$ lies above a real place of $K'$, and
\item\label{L2} $S=\{\rho\in\Sigma:\ \rho|_{K'}\in S'\}$.
\end{enumerate}
Without (L1), the set defined in (L2) would not be a system of representatives. We say that $S$ is
\emph{lifted from $K'$} if it is lifted from some $S'$; then $S'=\{\iota:\Sigma_\iota\subseteq S\}$ is
determined by $S$. Every $S$ is lifted from $K$ itself. In terms of blocks over $K'$, $S$ is lifted from $K'$ if
and only if every block lies in $S$ or in $\Sigma\setminus S$, and every block mapped to itself by $c$ consists
of real embeddings. If (L1) holds, the systems lifted from $K'$ are in bijection with the systems of
representatives of $K'$, so there are $2^{r_2(K')}$ of them.

\smallskip\noindent\emph{Admissible pairs.} A pair $(K',K_0)$ of subfields with $K_0\subseteq K'$ is
\emph{admissible} if $K_0$ is totally real and above each place of $K_0$ there is at most one complex place
of $K'$. Since $\Q$ is totally real with one place, $(K',\Q)$ is admissible exactly when $r_2(K')\le1$.

\smallskip
We keep to the following use of letters: $K'$ is a subfield from which $S$ is lifted, $K_0$ is the totally
real member of an admissible pair, and $K^{+}$ is always the largest totally real subfield of $K$. The letter
$K_c$ is reserved for one particular subfield, constructed from $G$ and $c$ in Section~\ref{sec:canonical}.

\begin{theorem}[the complete list]\label{thm:list}
Let $K$ be a number field with $r_2\ge1$ and $S$ a system of representatives. The following are equivalent.
\begin{enumerate}[(a),leftmargin=0.9cm,itemsep=2pt]
\item $T_{K,S}$ is quasi-abelian.
\item $S$ is lifted from a subfield $K'$ such that $(K',K^{+})$ is admissible. In particular
  $K^{+}\subseteq K'$, and above each place of $K^{+}$ there is at most one complex place of $K'$.
\item $S$ is lifted from a subfield $K'$ such that $(K',K_0)$ is admissible for some subfield
  $K_0\subseteq K'$.
\end{enumerate}
Moreover, suppose that some $T_{K,S}$ is quasi-abelian. Then (b) holds with one and the same subfield $K'$ for all
quasi-abelian $T_{K,S}$, and with no other subfield. We denote it by $K_c$; Section~\ref{sec:canonical} describes
it in terms of $G$ and $c$. The quasi-abelian $T_{K,S}$ are exactly those with $S$ lifted from $K_c$, and there are
$2^{r_2(K_c)}$ of them.
\end{theorem}

The implication (b)$\Rightarrow$(c) is clear, with $K_0=K^{+}$, and (c)$\Rightarrow$(a) combines the
trace-form criterion with Dodson's lifting theorem (Section~\ref{sec:subfields}). The substance is
(a)$\Rightarrow$(b), which takes up Sections~\ref{sec:graph}--\ref{sec:canonical}. Condition (c) contains the
two constructions of the introduction: $K'=K$ is the trace-form criterion, and $K_0=\Q$ means that $S$ is
lifted from a subfield $K'$ with $r_2(K')=1$. No assumption on $r_1$ is made. For $r_1=0$ the theorem gives
back Shimura--Taniyama (Section~\ref{sec:consequences}).

\begin{remark}[other lattices]\label{rem:lattices}
In both theorems the ring $\OK$ can be replaced by any finitely generated subgroup $\mathfrak a$ of $K$ of rank $\delta$, for
instance an order or a fractional ideal. Indeed, $\mathfrak a$ and $\OK$ are commensurable: each contains an integer multiple of the other. Hence the groups
$\C^{g}/u(\mathfrak a)$ and $T_{K,S}$ are isogenous, and a suitable integer multiple of an ample Riemann form
for one of them is an ample Riemann form for the other. These are the groups that occur in Dodson's result
quoted in the introduction~\cite{Dod80}. For fractional ideals they are also considered in~\cite{DF}.
\end{remark}

\section{Riemann forms as rational forms on \texorpdfstring{$K$}{K}}\label{sec:reduction}

We first move the problem from $\C^{g}$ to $K$. Let $F_K=u_\R^{-1}(F)\subseteq K\otimes\R$; it consists
of the $x$ with $\sigma_i(x)=0$ for all real $\sigma_i$. It carries the complex structure
$J=u_\R^{-1}\circ\iu\circ u_\R$, so that $\tau_j(Jx)=\iu\,\tau_j(x)$. For a $\Q$-bilinear form $E$ on $K$ we
write $E_\R$ for its $\R$-bilinear extension to $K\otimes\R$.

\begin{proposition}\label{prop:reduction}
$T_{K,S}$ is quasi-abelian if and only if there is an alternating $\Q$-bilinear form $E\colon K\times K\to\Q$
such that
\begin{description}[leftmargin=1.1cm,itemsep=2pt]
\item[(I)] $E(\OK\times\OK)\subseteq\Z$;
\item[(H)] $E_\R(Jx,Jy)=E_\R(x,y)$ for all $x,y\in F_K$;
\item[(P)] the form $H(x,y)=E_\R(Jx,y)+\iu\,E_\R(x,y)$ is positive definite on $F_K$.
\end{description}
\end{proposition}

Condition (H) says exactly that $H$ is Hermitian for the complex structure $J$. The proof rests on a small
extension lemma. We use coordinates $(\zeta,z)\in\C^{r_1}\times\C^{r_2}$ on $\C^{g}$, so that
$W:=\R\Lambda=\R^{r_1}\times\C^{r_2}$ and $F=\{0\}\times\C^{r_2}$. Hermitian forms are linear in the first
variable.

\begin{lemma}\label{lem:extension}
Let $\omega$ be a real alternating form on $W$ with $\omega(\iu v,\iu w)=\omega(v,w)$ for $v,w\in F$. Then
there is a Hermitian form $\widetilde H$ on $\C^{g}$ with $\im\widetilde H=\omega$ on $W\times W$ and
$\widetilde H(v,w)=\omega(\iu v,w)+\iu\,\omega(v,w)$ for $v,w\in F$.
\end{lemma}

\begin{proof}
We build $\widetilde H$ block by block. On $F$ put $H_F(z,z')=\omega(\iu z,z')+\iu\,\omega(z,z')$; the
invariance of $\omega$ under $\iu$ is exactly what makes $H_F$ Hermitian. On the real directions let
$A=(\omega(e_k,e_l))$, a real antisymmetric $r_1\times r_1$ matrix, and put
$H_1(\zeta,\zeta')=\iu\sum_{k,l}A_{kl}\,\zeta_k\bar\zeta'_l$; this is Hermitian, and
$\im H_1(s,s')=\omega(s,s')$ for real $s,s'$. For the mixed terms, $\lambda_k(z)=\omega(e_k,z)$ is a real
linear form on $\C^{r_2}$, and $\mu_k(z)=-\lambda_k(\iu z)+\iu\,\lambda_k(z)$ is the unique conjugate-linear
form with $\im\mu_k=\lambda_k$. Put $\beta(\zeta,z')=\sum_k\zeta_k\,\mu_k(z')$ and
\[
  \widetilde H\bigl((\zeta,z),(\zeta',z')\bigr)=H_1(\zeta,\zeta')+\beta(\zeta,z')+\overline{\beta(\zeta',z)}+H_F(z,z').
\]
This is Hermitian. For $s,s'$ real, $\im\beta(s,z')=\omega(s,z')$ and
$\im\overline{\beta(s',z)}=-\omega(s',z)=\omega(z,s')$,
so $\im\widetilde H=\omega$ on $W\times W$, and on $F$ only $H_F$ remains.
\end{proof}

\begin{proof}[Proof of Proposition~\ref{prop:reduction}]
Let $\widetilde H$ be an ample Riemann form and put $E(a,b)=\im\widetilde H(u(a),u(b))$ for $a,b\in K$. It
is alternating, integral on $\OK$, and therefore rational, since a multiple of any element of $K$ lies in
$\OK$. Its real extension is $E_\R(x,y)=\im\widetilde H(u_\R x,u_\R y)$. The identity
$\im\widetilde H(\iu v,\iu w)=\im\widetilde H(v,w)$ gives (H). From $\im\widetilde H(\iu v,w)=\re\widetilde H(v,w)$
we get $H(x,y)=\widetilde H(u_\R x,u_\R y)$ on $F_K$, which is positive: this is (P).

Conversely, let $E$ satisfy (I), (H) and (P). Transport $E_\R$ to $W$ by $u_\R$ and apply
Lemma~\ref{lem:extension}. The Hermitian form $\widetilde H$ obtained restricts on $F$ to the transport of
$H$, which is positive by (P), and $\im\widetilde H(\Lambda\times\Lambda)=E(\OK\times\OK)\subseteq\Z$ by (I).
\end{proof}

\section{Coordinates indexed by the embeddings}\label{sec:coordinates}

Fix an integral basis $\omega_1,\dots,\omega_\delta$ of $\OK$. It is a $\Q$-basis of $K$, and a basis of
$K\otimes_\Q k$ over $k$ for every field $k$ containing $\Q$. For $x\in K\otimes k$ we write $[x]\in k^{\delta}$
for its coordinate vector. Put
\[
  M=\bigl(\rho(\omega_a)\bigr)_{\rho\in\Sigma,\ 1\le a\le\delta}\in M_\delta(\Kt),
\]
with one row for each embedding. For $x\in K\otimes\C$ the vector $M[x]$ is $(\rho(x))_{\rho\in\Sigma}$,
where $\rho(x)$ denotes the $\C$-linear extension of $\rho$. The matrix $M$ is invertible, since
$(\det M)^{2}$ is the discriminant of $K$. Thus $x\mapsto M[x]$ identifies $K\otimes\C$ with $\C^{\Sigma}$. Let
$e_\rho\in K\otimes\Kt$ be the element with $M[e_\rho]$ the unit vector at $\rho$; its coordinate vector
$[e_\rho]$ is the column of $M^{-1}$ with index $\rho$. Every $x\in K\otimes\C$ is then
$x=\sum_\rho\rho(x)\,e_\rho$. (The $e_\rho$ are the primitive idempotents of the algebra $K\otimes\Kt$, but
we shall not need this.)

The Galois group acts on $M$ in the simplest possible way. For $\gamma\in G$, let $\gamma(M)$ be obtained by
applying $\gamma$ to every entry. Since $\gamma(\rho(\omega_a))=(\gamma\rho)(\omega_a)$, the row of $\gamma(M)$
with index $\rho$ is the row of $M$ with index $\gamma\rho$:
\begin{equation}\label{eq:permute}
  \gamma(M)=P_\gamma M ,
\end{equation}
where $P_\gamma$ is the permutation matrix with $(P_\gamma X)_{\rho}=X_{\gamma\rho}$ for the rows of any matrix
$X$ indexed by $\Sigma$. In words: \emph{$\gamma$ permutes the rows of $M$.} For $\gamma=c$ this reads
$\overline M=P_cM$.

\medskip\noindent\textbf{Entries of a form.} For a field $k\supseteq\Q$, an alternating $k$-bilinear form $E$ on
$K\otimes k$ is given by its Gram matrix $B=(E(\omega_a,\omega_b))\in\Alt_\delta(k)$, through
$E(x,y)=[x]^{\mathsf T}B\,[y]$; we write $E^{B}$ for the form with Gram matrix $B$. Rational forms on $K$
correspond to $\Alt_\delta(\Q)$ and real forms on $K\otimes\R$ to $\Alt_\delta(\R)$. Each of them extends to
$K\otimes\C$ with the same Gram matrix, and we use the same letter for the extension. The \emph{entry
matrix} of $E=E^{B}$ is
\[
  \mathcal E(B)=(M^{-1})^{\mathsf T}\,B\,M^{-1}=\bigl(E_{\rho\rho'}\bigr)_{\rho,\rho'\in\Sigma},\qquad
  E_{\rho\rho'}=E(e_\rho,e_{\rho'}),
\]
and since $[x]=M^{-1}(M[x])$,
\begin{equation}\label{eq:entries}
  E(x,y)=(M[x])^{\mathsf T}\,\mathcal E(B)\,(M[y])=\sum_{\rho,\rho'\in\Sigma}\rho(x)\,\rho'(y)\,E_{\rho\rho'} .
\end{equation}
Rationality is easy to see in $B$; the complex structure, as we shall see, is easy to see in
$\mathcal E(B)$. The entries are antisymmetric, $E_{\rho'\rho}=-E_{\rho\rho'}$, so it makes sense to say that
$E$ vanishes on an unordered pair. Explicitly, the definition of $\mathcal E(B)$ gives
\[
  E^{B}_{\rho\rho'}=\sum_{a,b=1}^{\delta}(M^{-1})_{a\rho}\,(M^{-1})_{b\rho'}\,B_{ab}.
\]
The entries $\rho(\omega_a)$ of $M$ lie in $\Kt$, hence so do those of $M^{-1}$. Thus each $E^{B}_{\rho\rho'}$ is
a linear combination of the entries $B_{ab}$ of the Gram matrix, with coefficients in $\Kt$ that do not depend
on $B$. In particular $E^{B}_{\rho\rho'}\in\Kt$ when $B$ is rational.

\medskip\noindent\textbf{The dictionary.} Let $x,y\in F_K$ and write $z_j=\tau_j(x)$, $w_k=\tau_k(y)$. Then
$\sigma_i(x)=0$ and $\bar\tau_j(x)=\bar z_j$, and $J$ multiplies $z_j$ by $\iu$ and $\bar z_j$ by $-\iu$. In
terms of $M$,
\[
  M[Jx]=\Delta_S\,M[x]\qquad(x\in F_K),
\]
where $\Delta_S$ is the diagonal matrix with entry $0$ at each $\sigma_i$, $\iu$ at each $\tau_j$ and $-\iu$ at
each $\bar\tau_j$.

\begin{lemma}\label{lem:dictionary}
Let $E$ be a real alternating form on $K\otimes\R$.
\begin{enumerate}[(i),leftmargin=0.9cm,itemsep=2pt]
\item $E_{c\rho,c\rho'}=\overline{E_{\rho\rho'}}$ for all $\rho,\rho'$.
\item For $x,y\in F_K$,
  \[
    E(Jx,Jy)-E(x,y)=-2\sum_{j,k}\Bigl(E_{\tau_j\tau_k}\,z_jw_k+E_{\bar\tau_j\bar\tau_k}\,\bar z_j\bar w_k\Bigr).
  \]
  Hence $E$ satisfies (H) if and only if $E$ vanishes on all forbidden pairs.
\item If $E$ satisfies (H), then $H(x,y)=E(Jx,y)+\iu E(x,y)=\sum_{j,k}h_{jk}\,z_j\bar w_k$ with
  $h_{jk}=2\iu\,E_{\tau_j\bar\tau_k}$. Hence (P) holds if and only if the Hermitian matrix $h=(h_{jk})$ is
  positive definite. In particular, (P) implies $E_{\tau_j\bar\tau_j}\ne0$ for every $j$: $E$ does not
  vanish on any diagonal pair.
\end{enumerate}
\end{lemma}

\begin{proof}
(i) The Gram matrix $B$ is real. Taking complex conjugates in $\mathcal E(B)=(M^{-1})^{\mathsf T}BM^{-1}$
and using $\overline M=P_cM$, hence $\overline{M^{-1}}=M^{-1}P_c^{\mathsf T}$, we get
\[
  \overline{\mathcal E(B)}=P_c\,(M^{-1})^{\mathsf T}BM^{-1}P_c^{\mathsf T}=P_c\,\mathcal E(B)\,P_c^{\mathsf T},
\]
which is (i) entry by entry.

(ii) By~\eqref{eq:entries} and $M[Jx]=\Delta_SM[x]$,
\[
  E(Jx,Jy)-E(x,y)=(M[x])^{\mathsf T}\bigl(\Delta_S\,\mathcal E\,\Delta_S-\mathcal E\bigr)(M[y]).
\]
The coordinates of $x$ and $y$ at the real embeddings vanish. At $(\tau_j,\bar\tau_k)$ and $(\bar\tau_j,\tau_k)$
the matrix $\Delta_S\mathcal E\Delta_S$ has the same entries as $\mathcal E$, since $\iu\cdot(-\iu)=1$; at
$(\tau_j,\tau_k)$ and $(\bar\tau_j,\bar\tau_k)$ the signs change. This gives the formula. Its right-hand side
vanishes for all $z,w$ if and only if all $E_{\tau_j\tau_k}$ and $E_{\bar\tau_j\bar\tau_k}$ vanish, since the
two sums are holomorphic and antiholomorphic in $(z,w)$. By (i) the latter are the conjugates of the former,
and $E_{\tau_j\tau_j}=0$ anyway. So (H) means $E_{\tau_j\tau_k}=0$ for $j\ne k$.

(iii) In the same way
\[
  H(x,y)=E(Jx,y)+\iu E(x,y)=(M[x])^{\mathsf T}\,(\Delta_S+\iu I)\,\mathcal E\,(M[y]).
\]
The diagonal matrix $\Delta_S+\iu I$ is $2\iu$ at the $\tau_j$ and $0$ at the $\bar\tau_j$, and the coordinates
of $x$ at the real embeddings vanish. So only the rows of the $\tau_j$ contribute: $H(x,y)=2\iu\sum_{j,k}\bigl(E_{\tau_j\tau_k}z_jw_k+E_{\tau_j\bar\tau_k}z_j\bar w_k\bigr)$.
Under (H) the first sum is zero.
\end{proof}

To summarize: condition (H) asks for zeros on the forbidden pairs, and condition (P) needs nonzero values
on the diagonal pairs. The question is whether a \emph{rational} form can do both. The next lemma is the only
place where rationality enters.

\begin{lemma}[Galois equivariance]\label{lem:equivariance}
Let $B\in\Alt_\delta(\Kt)$ and $\gamma\in G$, and let $\gamma(B)$ be obtained by applying $\gamma$ to every
entry. Then
\[
  E^{\gamma(B)}_{\gamma\rho,\gamma\rho'}=\gamma\bigl(E^{B}_{\rho\rho'}\bigr)\qquad\text{for all }\rho,\rho'\in\Sigma.
\]
In particular, if $E$ is a rational form, then $E_{\gamma\rho,\gamma\rho'}=\gamma(E_{\rho\rho'})$, and if $E$
vanishes on one pair, it vanishes on the whole $G$-orbit of that pair.
\end{lemma}

\begin{proof}
Apply $\gamma$ to every entry of the matrices involved; this commutes with sums, products, transposes and
inverses. By~\eqref{eq:permute}, $\gamma(M^{-1})=(P_\gamma M)^{-1}=M^{-1}P_\gamma^{\mathsf T}$, and therefore
\[
  \gamma\bigl(\mathcal E(B)\bigr)=\bigl(M^{-1}P_\gamma^{\mathsf T}\bigr)^{\mathsf T}\,\gamma(B)\,M^{-1}P_\gamma^{\mathsf T}
  =P_\gamma\,\mathcal E\bigl(\gamma(B)\bigr)\,P_\gamma^{\mathsf T}.
\]
The entry with index $(\rho,\rho')$ of the right-hand side is $E^{\gamma(B)}_{\gamma\rho,\gamma\rho'}$. If $B$ is
rational, then $\gamma(B)=B$.
\end{proof}

\section{Proof of the orbit criterion}\label{sec:proof}

\subsection*{(b) \texorpdfstring{$\Leftrightarrow$}{<=>} (c)}
The $G$-orbits of the diagonal pairs consist of the pairs
$\gamma\{\tau,c\tau\}=\{\gamma\tau,\,(\gamma c\gamma^{-1})\gamma\tau\}$, where $\tau$ is a complex embedding.
These are exactly the pairs $\{\rho,c'\rho\}$ with $c'$ a complex conjugation of $\Kt$ and $c'\rho\ne\rho$.
Hence a forbidden pair $\{\tau_j,\tau_k\}$ lies in the orbit of a diagonal pair if and only if some complex
conjugation maps $\tau_j$ to $\tau_k$.

\subsection*{(a) \texorpdfstring{$\Rightarrow$}{=>} (b)}
Let $T_{K,S}$ be quasi-abelian, and let $E$ be the rational form given by
Proposition~\ref{prop:reduction}. Suppose that some $\gamma\in G$ maps a forbidden pair
$\{\tau_k,\tau_l\}$ onto a diagonal pair $\{\tau_j,\bar\tau_j\}$. By Lemma~\ref{lem:dictionary}(ii), $E$
vanishes on $\{\tau_k,\tau_l\}$. By Lemma~\ref{lem:equivariance} it then vanishes on $\{\tau_j,\bar\tau_j\}$,
which contradicts Lemma~\ref{lem:dictionary}(iii).

\subsection*{(b) \texorpdfstring{$\Rightarrow$}{=>} (a)}
This direction needs a construction. Let $\mathcal P$ be the union of the $G$-orbits of pairs that contain
no forbidden pair. By (b), every diagonal pair lies in $\mathcal P$. Consider the spaces of forms that vanish
outside $\mathcal P$:
\begin{align*}
  V&=\{B\in\Alt_\delta(\Q):\ E^{B}_{\rho\rho'}=0\ \text{for all pairs }\{\rho,\rho'\}\notin\mathcal P\},\\
  Y&=\{B\in\Alt_\delta(\R):\ E^{B}_{\rho\rho'}=0\ \text{for all pairs }\{\rho,\rho'\}\notin\mathcal P\}.
\end{align*}
Every forbidden pair lies outside $\mathcal P$, so by Lemma~\ref{lem:dictionary}(ii) every form in $V$ or
in $Y$ satisfies (H). We show that $V$ is dense in $Y$, and that $Y$ contains a form satisfying (P). A
rational form close to it will do.

We use a standard descent lemma.

\begin{lemma}[Galois descent]\label{lem:descent}
Let $L/\Q$ be a finite Galois extension with group $\Gamma_L$, and let $U\subseteq L^{m}$ be an $L$-subspace
with $\gamma(U)=U$ for all $\gamma\in\Gamma_L$, where $\gamma$ acts on coordinates. Then $U$ is spanned over $L$
by $U\cap\Q^{m}$.
\end{lemma}

\begin{proof}
We apply the trace $\Tr_{L/\Q}$ to vectors in $L^{m}$ coordinatewise. Since $L/\Q$ is Galois,
$\Tr_{L/\Q}(\xi)=\sum_{\gamma\in\Gamma_L}\gamma(\xi)$ for $\xi\in L^{m}$, and this vector lies in $\Q^{m}$.

Let $(\kappa_i)$ be a basis of $L$ over $\Q$. The trace form of $L$ over $\Q$ is nondegenerate, so there is a
dual basis $(\kappa_i^{*})$ of $L$ over $\Q$, with $\Tr_{L/\Q}(\kappa_i\kappa_j^{*})=1$ if $i=j$ and $0$
otherwise. Every $\xi\in L^{m}$ can be written $\xi=\sum_j\kappa_j^{*}q_j$ with $q_j\in\Q^{m}$, and then
$\Tr_{L/\Q}(\kappa_i\xi)=\sum_j\Tr_{L/\Q}(\kappa_i\kappa_j^{*})\,q_j=q_i$. Hence
\[
  \xi=\sum_i\kappa_i^{*}\,\Tr_{L/\Q}(\kappa_i\xi)\qquad(\xi\in L^{m}).
\]
Now let $\xi\in U$. For each $i$ the vector $\kappa_i\xi$ lies in $U$, hence so does $\gamma(\kappa_i\xi)$ for
every $\gamma\in\Gamma_L$, and therefore $\Tr_{L/\Q}(\kappa_i\xi)\in U\cap\Q^{m}$. By the formula above, $\xi$
is an $L$-linear combination of these vectors. Since $U\cap\Q^{m}$ is contained in $U$, its $L$-span is
exactly $U$.
\end{proof}

\noindent\emph{Step 1: $V$ is dense in $Y$.} Let $U\subseteq\Alt_\delta(\Kt)$ be defined by the same
equations as $V$ and $Y$. The set of pairs outside $\mathcal P$ is $G$-stable, so Lemma~\ref{lem:equivariance}
shows that $\gamma(U)=U$ for all $\gamma\in G$. By Lemma~\ref{lem:descent}, $U$ is spanned over $\Kt$ by
$U\cap\Alt_\delta(\Q)=V$. The equations have coefficients in $\Kt$, and the rank of a linear system does not change when the field is
enlarged. So their solution space in $\Alt_\delta(\C)$ is the $\C$-span of $U$, hence of $V$. Now take a $\Q$-basis $(v_i)$ of $V$; it
stays linearly independent over $\C$. If a real matrix $B=\sum_it_iv_i$, with $t_i\in\C$, lies in $Y$, then
$\sum_i(t_i-\bar t_i)v_i=B-\bar B=0$, so all $t_i$ are real. Hence $Y$ is the $\R$-span of $V$, and $V$ is
dense in $Y$.

\medskip\noindent\emph{Step 2: a positive real form.} On $K\otimes\R$ put
\[
  E^{0}(x,y)=\sum_{j=1}^{r_2}\im\bigl(\tau_j(x)\,\overline{\tau_j(y)}\bigr)
  =\frac{1}{2\iu}\sum_{j=1}^{r_2}\bigl(\tau_j(x)\bar\tau_j(y)-\bar\tau_j(x)\tau_j(y)\bigr).
\]
The right-hand side is already $\C$-bilinear in the coordinates, so by~\eqref{eq:entries} the entries of
$E^{0}$ are $E^{0}_{\tau_j\bar\tau_j}=\frac{1}{2\iu}$, $E^{0}_{\bar\tau_j\tau_j}=-\frac{1}{2\iu}$, and $0$
elsewhere. Thus $E^{0}$ is supported on the diagonal pairs, which lie in $\mathcal P$: its Gram matrix lies
in $Y$. Its Hermitian matrix is $h=(2\iu E^{0}_{\tau_j\bar\tau_k})=I$, that is,
$H(x,y)=\sum_jz_j\bar w_j$.

\medskip\noindent\emph{Step 3: conclusion.} By Lemma~\ref{lem:dictionary}(iii) the map
$B\mapsto h(B)=(2\iu E^{B}_{\tau_j\bar\tau_k})$ is an $\R$-linear map from $Y$ to the Hermitian
$r_2\times r_2$ matrices, and $h=I$ at the Gram matrix of $E^{0}$. Positive definite matrices form an open
set. By Step~1 there is a rational $B\in V$ so close to the Gram matrix of $E^{0}$ that $h(B)$ is still
positive definite. The form $E^{B}$ is rational and satisfies (H) and (P). The entries of $B$ have a
common denominator $d$, and since the $\omega_a$ form a $\Z$-basis of $\OK$, the form $dE^{B}$ is integral
on $\OK$. Thus $dE^{B}$ satisfies (I), (H) and (P), and $T_{K,S}$ is quasi-abelian by
Proposition~\ref{prop:reduction}. \qed

\section{The N\'eron--Severi group}\label{sec:NS}

Let $X=\C^{n}/\Lambda_X$ be a toroidal group. The N\'eron--Severi group $\NS(X)$ of $X$ is the group of first Chern classes of holomorphic line bundles on $X$. It can be identified with the group of alternating
forms $E\colon\Lambda_X\times\Lambda_X\to\Z$ whose real extension to $\R\Lambda_X$ satisfies
$E(\iu v,\iu w)=E(v,w)$ for $v,w$ in the largest complex subspace of $\R\Lambda_X$;
see~\cite{AK} and~\cite{Um03,Um17}. An ample Riemann form gives an element of $\NS(X)$, but in general an element of
$\NS(X)$ need not be positive anywhere.

For $X=T_{K,S}$ we transport everything to $K$ by $u$, as in Section~\ref{sec:reduction}. Then the group
$\NS(T_{K,S})$ consists of the alternating forms on $\OK$ with values in $\Z$ that satisfy (H). Since every
rational form has a multiple that is integral on $\OK$, the space $\NS(T_{K,S})\otimes\Q$ consists of the
rational alternating forms on $K$ with (H). The proof of Theorem~\ref{thm:main} computes this space.

\begin{proposition}\label{prop:NS}
For an orbit $O$ of pairs, let $V_O$ be the space of rational alternating forms on $K$ whose entries vanish
outside $O$. Then $\dim_\Q V_O=|O|$, and
\[
  \NS(T_{K,S})\otimes\Q=\bigoplus_{O}V_O ,
\]
where $O$ runs through the $G$-orbits of pairs that contain no forbidden pair. In particular
\[
  \operatorname{rank}\NS(T_{K,S})=\sum_{O\ \text{without forbidden pairs}}|O| .
\]
\end{proposition}

\begin{proof}
By Lemmas~\ref{lem:dictionary}(ii) and~\ref{lem:equivariance}, a rational form satisfies (H) if and only if it
vanishes on every orbit that contains a forbidden pair. Let $E=E^{B}$ be such a form and $O$ an orbit. Let
$\mathcal E_O$ be the matrix that agrees with $\mathcal E(B)$ at the ordered pairs $(\rho,\rho')$ with
$\{\rho,\rho'\}\in O$ and is $0$ elsewhere, and put $B_O=M^{\mathsf T}\mathcal E_OM$. Since $O$ is $G$-stable,
Lemma~\ref{lem:equivariance} gives $\gamma(\mathcal E_O)=P_\gamma\mathcal E_OP_\gamma^{\mathsf T}$, and with
\eqref{eq:permute}
\[
  \gamma(B_O)=M^{\mathsf T}P_\gamma^{\mathsf T}\,P_\gamma\mathcal E_OP_\gamma^{\mathsf T}\,P_\gamma M=B_O
  \qquad(\gamma\in G).
\]
Hence $B_O$ is rational, and $B=\sum_OB_O$. The sum is direct, because forms supported on different orbits
have disjoint sets of nonzero entries. This proves the decomposition.

For the dimension, let $U_O\subseteq\Alt_\delta(\Kt)$ be the space of matrices $B$ with $\mathcal E(B)$
supported on $O$. Since $B\mapsto\mathcal E(B)$ is a linear bijection onto the alternating matrices indexed
by $\Sigma$, and such a matrix supported on $O$ has one free entry for each pair in $O$, we have
$\dim_{\Kt}U_O=|O|$. By Lemma~\ref{lem:equivariance}, $U_O$ is stable under $G$, so Lemma~\ref{lem:descent}
shows that it is spanned by $V_O=U_O\cap\Alt_\delta(\Q)$. Rational matrices that are linearly independent over
$\Q$ remain so over $\Kt$; hence $\dim_\Q V_O=\dim_{\Kt}U_O=|O|$.
\end{proof}

In particular, $\NS(T_{K,S})\ne0$ as soon as one orbit contains no forbidden pair, while $T_{K,S}$ is quasi-abelian
only if these orbits contain all the diagonal pairs. Examples where the rank is positive but $T_{K,S}$ is not
quasi-abelian are given in Section~\ref{sec:examples}.

\section{Totally real subfields and lifting}\label{sec:subfields}

Let $K_0\subseteq K$ be totally real, with embeddings $\iota_1,\dots,\iota_m$. Every block
$\Sigma_i=\Sigma_{\iota_i}$ over $K_0$ is mapped to itself by $c$, so it consists of $a_i$ real embeddings and
$b_i$ pairs of complex conjugate ones, where $a_i+2b_i=[K:K_0]$. The number $b_i$ counts the complex places of
$K$ above the $i$-th place of $K_0$. Thus $(K,K_0)$ is admissible exactly when $b_i\le1$ for all $i$.

\begin{theorem}[trace-form criterion]\label{thm:trace}
If $K$ contains a totally real subfield $K_0$ with $b_i\le1$ for all $i$, that is, if $(K,K_0)$ is
admissible, then $T_{K,S}$ is quasi-abelian for every system of representatives $S$.
\end{theorem}

\begin{proof}
Call a pair of embeddings \emph{inside} if both lie in the same block over $K_0$, and \emph{across} otherwise.
Since $G$ permutes the blocks, every $G$-orbit of pairs consists either of inside pairs only or of across pairs
only. A diagonal pair $\{\tau_j,c\tau_j\}$ is inside, because $c$ maps every block to itself. A forbidden pair
$\{\tau_j,\tau_k\}$ is across: if $\tau_j$ and $\tau_k$ lay in the same block $\Sigma_i$, then so would
$\bar\tau_j$ and $\bar\tau_k$, and $b_i\ge2$. Hence condition (b) of Theorem~\ref{thm:main} holds.
\end{proof}

The theorem contains the families found earlier. There are three of them. Dodson treated quadratic
extensions of a totally real field~\cite{Dod79}; there $b_i\le1$ is automatic. Andreotti and Gherardelli
treated cubic extensions $K/K_0$ that are not totally real at any place of $K_0$: above each place of $K_0$
lie one real and one complex place of $K$~\cite[\S4.2]{AG73}; see also~\cite{AG}. Abe treated
extensions of any relative degree with exactly one complex place above each place of $K_0$~\cite{Abe13}. In
the last two families $b_i=1$ for all $i$: the extension is uniform over $K_0$. Theorem~\ref{thm:trace} asks
only for $b_i\le1$, so $b_i=0$ may occur at some places of $K_0$ and $b_i=1$ at others. For $K_0=\Q$ the theorem reduces to
Corollary~\ref{cor:easy}(i). These results were originally proved with Riemann forms of the shape
$\Tr_{K_0/\Q}\circ\Phi$, which give the criterion its name. The proof above makes the choice of a form
unnecessary.

\medskip
We turn to lifting, defined in Section~\ref{subsec:subfields}. If $S$ is lifted from a system of
representatives $S'$ of $K'$, then the complex places of $K$ lie above complex places of $K'$ by (L1); so
$r_2(K')\ge1$, and $T_{K',S'}$ is defined.

\begin{theorem}[Dodson~{\cite[Thm.~3.1]{Dod79}}]\label{thm:lifting}
Let $S$ be lifted from a system of representatives $S'$ of a subfield $K'$, and let $m=[K:K']$.
\begin{enumerate}[(a),leftmargin=0.9cm,itemsep=2pt]
\item $T_{K,S}$ is isogenous to $(T_{K',S'})^{m}$; it is isomorphic to it if $\OK$ is a free
  $\mathcal O_{K'}$-module.
\item If $T_{K',S'}$ is quasi-abelian, then so is $T_{K,S}$.
\end{enumerate}
\end{theorem}

\begin{proof}
(a) Choose $\beta_1,\dots,\beta_m\in\OK$ forming a basis of $K$ over $K'$, and an $\mathcal O_{K'}$-basis of
$\OK$ if there is one. For $\iota\in S'$ let $\rho_{\iota,1},\dots,\rho_{\iota,m}$ be the embeddings of $K$
that extend $\iota$. By (L2) they lie in $S$, and as $\iota$ runs through $S'$ they make up all of $S$. Order
the coordinates of $\C^{g}$ accordingly. For $x=\sum_i a_i\beta_i$ with $a_i\in K'$,
\[
  \rho_{\iota,k}(x)=\sum_{i=1}^{m}\rho_{\iota,k}(\beta_i)\,\iota(a_i).
\]
Let $A\colon(\C^{S'})^{m}\to\C^{g}$ be the $\C$-linear map that sends $(v_1,\dots,v_m)$ to the vector with
coordinates $\sum_i\rho_{\iota,k}(\beta_i)\,(v_i)_\iota$. Then
\[
  A\bigl(u'(a_1),\dots,u'(a_m)\bigr)=u\Bigl(\sum_i a_i\beta_i\Bigr),
\]
where $u'$ is the map $u$ for $(K',S')$. The map $A$ is block diagonal, with one $m\times m$ block
$B_\iota=(\rho_{\iota,k}(\beta_i))_{k,i}$ for each $\iota\in S'$, and
$(\det B_\iota)^{2}=\iota\bigl(\operatorname{disc}_{K/K'}(\beta_1,\dots,\beta_m)\bigr)\ne0$. Hence $A$ is an
isomorphism. It maps $\Lambda'^{\,m}$, where $\Lambda'=u'(\mathcal O_{K'})$, onto
$\Lambda_1=u\bigl(\bigoplus_i\mathcal O_{K'}\beta_i\bigr)$, a subgroup of finite index in $\Lambda$. Hence $A$
induces an isomorphism $(T_{K',S'})^{m}\cong\C^{g}/\Lambda_1$, followed by the isogeny
$\C^{g}/\Lambda_1\to\C^{g}/\Lambda$ with kernel $\Lambda/\Lambda_1$. If the $\beta_i$ form an
$\mathcal O_{K'}$-basis, then $\Lambda_1=\Lambda$.

(b) Let $\widetilde H'$ be an ample Riemann form for $T_{K',S'}$. The orthogonal sum of $m$ copies is an ample
Riemann form for $(T_{K',S'})^{m}$. Transport it by $A$ to a Hermitian form $\widetilde H_1$ on $\C^{g}$. Then
$\widetilde H_1$ is positive on the largest complex subspace of $\R\Lambda_1=\R\Lambda$, which is $F$, and
$\im\widetilde H_1(\Lambda_1\times\Lambda_1)\subseteq\Z$. If $e$ is the index of $\Lambda_1$ in $\Lambda$, then
$e\Lambda\subseteq\Lambda_1$. Since
$\im\bigl(e^{2}\widetilde H_1\bigr)(\lambda,\mu)=\im\widetilde H_1(e\lambda,e\mu)\in\Z$ for $\lambda,\mu\in\Lambda$,
the form $e^{2}\widetilde H_1$ is an ample Riemann form for $T_{K,S}$.
\end{proof}

Part (b) can also be read off from condition (c) of the orbit criterion, which shows how this criterion
handles lifting.

\begin{proof}[Second proof of (b)]
Suppose a complex conjugation $c'=\gamma c\gamma^{-1}$ interchanges $\tau_j\ne\tau_k$ in $S$, and let
$\tau_j\in\Sigma_\iota$. If $c'$ maps the block $\Sigma_\iota$ to itself, then $c$ fixes $\gamma^{-1}\iota$, which
is therefore a real embedding of $K'$. By (L1) all elements of $\Sigma_{\gamma^{-1}\iota}$ are real, i.e.\ fixed
by $c$; as $\Sigma_\iota=\gamma\Sigma_{\gamma^{-1}\iota}$, the conjugate $c'$ fixes $\Sigma_\iota$ pointwise, so
$c'\tau_j=\tau_j$, a contradiction. Otherwise $\tau_k$ lies in the block $\Sigma_{c'\iota}\ne\Sigma_\iota$, and
$\iota,c'\iota\in S'$ by (L2). Both are complex, by (L1), since $\tau_j$ and $\tau_k$ are. On the Galois
closure of $K'$, the restriction of $c'$ is a complex conjugation interchanging $\iota$ and $c'\iota$, so
condition (c) does not hold for $(K',S')$, and $T_{K',S'}$ is not quasi-abelian by Theorem~\ref{thm:main}.
\end{proof}

\begin{proof}[Proof of Theorem~\ref{thm:list}, (c)$\Rightarrow$(a)]
Let $(K',K_0)$ be admissible and $S$ lifted from $S'$. Theorem~\ref{thm:trace}, applied to the field $K'$ and
its totally real subfield $K_0$, shows that $T_{K',S'}$ is quasi-abelian; its hypothesis is exactly the
admissibility of $(K',K_0)$. By Theorem~\ref{thm:lifting}(b), $T_{K,S}$ is quasi-abelian.
\end{proof}

\section{The conjugation graph}\label{sec:graph}

We now prove the implication (a)$\Rightarrow$(b) of Theorem~\ref{thm:list}. Its input is the orbit criterion
in the form (c) of Theorem~\ref{thm:main}, together with Lemma~\ref{lem:N}. Riemann forms play no further
role: the rest is a study of the complex conjugations as permutations of $\Sigma$.

We use only a few notions from graph theory, all for graphs without loops or multiple edges. Two vertices are
\emph{adjacent}, or \emph{neighbours}, if an edge joins them. For a set $X$ of vertices, the \emph{induced
subgraph} on $X$ has vertex set $X$ and contains every edge of the graph with both ends in $X$. A
\emph{walk} is a finite sequence of vertices in which any two consecutive ones are adjacent, and a \emph{path}
is a walk without repeated vertices. A graph is \emph{connected} if any two of its vertices are joined by a
walk. Its \emph{connected components} are the maximal connected induced subgraphs; we use the same word for
their vertex sets. A \emph{triangle} is a set of three pairwise adjacent vertices. An \emph{automorphism} of a
graph is a bijection of its vertex set that maps edges to edges.

A \emph{proper $2$-colouring} of a graph is a partition of its vertex set into two sets, the \emph{colour
classes}, such that no edge joins two vertices of the same class. A graph is \emph{bipartite} if it has a
proper $2$-colouring. A bipartite graph contains no triangle, since two of the three vertices would have the
same colour. A connected bipartite graph has only one proper $2$-colouring, up to exchanging the two classes:
the class of one vertex determines the classes of its neighbours, then those of their neighbours, and so on
along walks.

Recall from Section~\ref{subsec:galois} that $\cD=\{\gamma c\gamma^{-1}:\gamma\in G\}$ is the set of complex
conjugations of $\Kt$, a conjugacy class of involutions of $G$, and that $N=\langle\cD\rangle$ is the subgroup it
generates. Since $\cD$ is stable under conjugation, $N$ is normal in $G$; its orbits on $\Sigma$ are the blocks over
$K^{+}$ (Lemma~\ref{lem:N}).

\begin{definition}
The \emph{conjugation graph} $\Gamma$ has vertex set $\Sigma$; two distinct embeddings are adjacent if some
$d\in\cD$ interchanges them. We write $\rho\sim\rho'$. For $d\in\cD$ let $R(d)$ be the set of embeddings fixed by
$d$ and $C(d)$ the set of embeddings moved by $d$. We put $R=R(c)$ and $C=C(c)$: the real and the complex
embeddings. For $X\subseteq\Sigma$, $\Gamma[X]$ denotes the induced subgraph.
\end{definition}

The vertices of $\Gamma$ are the embeddings, and we keep this word throughout. Thus a \emph{real} or
\emph{complex} embedding is one that lies in $R$ or in $C$, that is, one fixed or moved by $c$; for another
$d\in\cD$ we say explicitly that an embedding is fixed or moved by $d$. Adjacency and neighbours always refer to
$\Gamma$.

Since $\cD$ is a conjugacy class, every $\gamma\in G$ is an automorphism of $\Gamma$. In particular $c$ is one.
Each pair $\{\tau,c\tau\}$ with $\tau\in C$ is an edge.

\begin{lemma}\label{lem:components}
\begin{enumerate}[(a),leftmargin=0.9cm,itemsep=2pt]
\item Every embedding of $\Sigma$ is moved by some $d\in\cD$. In particular every $N$-orbit has at least two
  embeddings.
\item The connected components of $\Gamma$ are the $N$-orbits, that is, the blocks of $\Sigma$ over $K^{+}$.
\end{enumerate}
\end{lemma}

\begin{proof}
(a) An embedding fixed by every element of $\cD$ is fixed by $N$. Since $N$ is normal, the set of embeddings fixed by
$N$ is mapped to itself by $G$, so it is empty or all of $\Sigma$, as $G$ is transitive. It is not all of
$\Sigma$, because $c\in N$ moves the complex embeddings, and $r_2\ge1$.

(b) An edge joins $\rho$ to $d\rho$ with $d\in N$, so it stays inside an $N$-orbit. Conversely, every element of
$N$ is a product $d_k\cdots d_1$ of elements of $\cD$, and the sequence $\rho,\ d_1\rho,\ d_2d_1\rho,\dots$ becomes a walk in
$\Gamma$ once repeated consecutive terms are removed; they occur where some $d_i$ fixes the current
embedding. The last assertion is
Lemma~\ref{lem:N}.
\end{proof}

The orbit criterion becomes a statement about colourings.

\begin{proposition}\label{prop:colouring}
$T_{K,S}$ is quasi-abelian if and only if $C=(S\cap C)\sqcup(cS\cap C)$ is a proper $2$-colouring of
$\Gamma[C]$. Hence some $T_{K,S}$ is quasi-abelian if and only if $\Gamma[C]$ is bipartite, and then exactly
$2^{k}$ of the $2^{r_2}$ systems of representatives $S$ give a quasi-abelian $T_{K,S}$, where $k$ is the number of
connected components of $\Gamma[C]$.
\end{proposition}

\begin{proof}
By Theorem~\ref{thm:main}(c), $T_{K,S}$ is quasi-abelian if and only if no edge of $\Gamma$ joins two embeddings of
$S\cap C$. Applying the automorphism $c$, the same then holds for $cS\cap C=c(S\cap C)$, and these two sets
partition $C$. Conversely, if $C=X\sqcup Y$ is a proper $2$-colouring, then each edge $\{\tau,c\tau\}$ has one
end in $X$, so $S=R\cup X$ is a system of representatives, and condition (c) holds for it. A proper
$2$-colouring of a bipartite graph is determined by the colour of one embedding in each component, as
explained above; this gives the number $2^{k}$.
\end{proof}

From now on we assume
\[
  (\star)\qquad \Gamma[C]\ \text{is bipartite.}
\]
If $d=\gamma c\gamma^{-1}$, then $C(d)=\gamma C$ and $\Gamma[C(d)]=\gamma\,\Gamma[C]$. Hence $(\star)$ holds for every
$d\in\cD$ in place of $c$, and all statements below that are proved for $c$ hold for every $d\in\cD$. We often
use $(\star)$ in the following weak form: \emph{three pairwise adjacent embeddings are never all moved by the same
complex conjugation, say $d$}, since they would form a triangle in $\Gamma[C(d)]$.

\subsection*{Triangles}

The assumption $(\star)$ concerns only the induced subgraphs $\Gamma[C(d)]$, $d\in\cD$. The graph $\Gamma$ itself
may well contain triangles, and the next lemma produces them: as soon as a complex conjugation interchanges an
embedding moved by $d$ with one fixed by $d$, we obtain a triangle in $\Gamma$ that passes through an embedding
fixed by $d$. Thus $\Gamma$ is in general not bipartite, even under $(\star)$.

\begin{lemma}[triangle lemma]\label{lem:triangle}
Let $d,d'\in\cD$, and suppose that $d'$ interchanges an embedding $x\in C(d)$ with an embedding $r\in R(d)$. Then $r$,
$x$ and $dx$ are pairwise adjacent, and $d'$ fixes $dx$.
\end{lemma}

\begin{proof}
$d'$ interchanges $r$ and $x$, and $d$ interchanges $x$ and $dx$. The conjugate $dd'd$ sends $dx$ to
$dd'x=dr=r$, so it interchanges $r$ and $dx$. Thus $\{r,x,dx\}$ is a triangle in $\Gamma$; it is not contained in
$\Gamma[C(d)]$, since $d$ fixes $r$. The conjugation $d'$ moves $r$ and $x$; if it also moved $dx$, the triangle
would lie in $\Gamma[C(d')]$, which $(\star)$ for $d'$ excludes.
\end{proof}

The next lemma says that a complex conjugation which moves a complex embedding to a complex embedding leaves the
nearby real embeddings fixed.

\begin{lemma}\label{lem:fix}
Let $d\in\cD$ and $x\in C$ with $dx\in C$ and $dx\ne x$. Then $d$ moves $cx$, and $d$ fixes every real embedding
adjacent to $x$.
\end{lemma}

\begin{proof}
Suppose $d$ fixes $cx$. Then $c$ interchanges $x\in C(d)$ and $cx\in R(d)$, and Lemma~\ref{lem:triangle}, applied
to the pair $(d,c)$, says that $c$ fixes $dx$. But $dx$ is complex. So $d$ moves both $x$ and $cx$.

Now let $r\in R$ be adjacent to $x$, and let $d_1\in\cD$ interchange them. By Lemma~\ref{lem:triangle}, applied to
$(c,d_1)$, the embeddings $r$, $x$, $cx$ are pairwise adjacent. If $d$ moved $r$, it would move all three.
\end{proof}

\section{The structure of an \texorpdfstring{$N$}{N}-orbit}\label{sec:orbits}

We keep the assumption $(\star)$.

\begin{lemma}\label{lem:connected}
Let $O$ be an $N$-orbit with $C\cap O\ne\emptyset$. Then $\Gamma[C\cap O]$ is connected, and every real embedding of
$O$ has a complex neighbour.
\end{lemma}

\begin{proof}
Fix $x\in C\cap O$, let $Z$ be its component in $\Gamma[C]$, and let $\hat Z$ be the union of $Z$ and of the real
embeddings that have a neighbour in $Z$. Then $\hat Z\subseteq O$ by Lemma~\ref{lem:components}, and every
neighbour of an embedding of $Z$ lies in $\hat Z$, since a complex one lies in $Z$. We show that
$d\hat Z\subseteq\hat Z$ for every $d\in\cD$. Then $\hat Z$ is stable under $N$, so $\hat Z=O$. Hence
$C\cap O=Z$, and every real embedding of $O$ has a neighbour in $Z$.

For $\rho\in Z$, the embedding $d\rho$ is $\rho$ or a neighbour of $\rho$, so it lies in $\hat Z$. Now let
$r\in\hat Z$ be real, let $a\in Z$ be a neighbour of $r$, and let $d_1\in\cD$ interchange $a$ and $r$. It suffices
to find a neighbour of $dr$ in $Z$. If $da$ is complex, then $da\in Z$, since it is $a$ or a complex neighbour of
$a$; and the conjugate $dd_1d$ interchanges $da$ and $dr$. If $da$ is real, Lemma~\ref{lem:triangle}, applied to
the pair $(c,d)$, shows that $d$ fixes $ca$. The conjugate $cd_1c$ interchanges $ca$ and $cr=r$, so $d(cd_1c)d$
interchanges $ca$ and $dr$; and $ca\in Z$, since it is a complex neighbour of $a$.
\end{proof}

The hypothesis $C\cap O\ne\emptyset$ is not always satisfied: an $N$-orbit may consist of real embeddings only,
as in Example~\ref{ex:deg9}. Such an orbit is still moved by other complex conjugations
(Lemma~\ref{lem:components}), and Lemma~\ref{lem:connected}, applied to them, covers it in
Section~\ref{sec:canonical}.

Consequently, for every $N$-orbit $O$ meeting $C$ the graph $\Gamma[C\cap O]$ has exactly one proper
$2$-colouring, up to exchanging the colours. We write $C\cap O=P_O\sqcup Q_O$ for the two colour classes. The
automorphism $c$ maps each embedding of $C\cap O$ to a neighbour, so it exchanges the two classes: $Q_O=cP_O$. Let $S$ be a system of representatives with $T_{K,S}$ quasi-abelian. By
Proposition~\ref{prop:colouring}, $S\cap C$ and $cS\cap C$ form a proper $2$-colouring of $\Gamma[C]$. By
Lemma~\ref{lem:connected}, $\Gamma[C\cap O]$ is connected, so it is a connected component of $\Gamma[C]$. Hence
$S\cap C\cap O$ is $P_O$ or $Q_O$: $S$ contains exactly one of the two classes.

For $\rho$ in an $N$-orbit $O$ put
\[
  \Pi(\rho)=\{\rho\}\cup\{\rho'\in O:\ \rho'\not\sim\rho\},
\]
the embedding together with its non-neighbours in its $N$-orbit. This set is defined by $\Gamma$ alone.

\begin{lemma}\label{lem:parts}
Let $O$ be an $N$-orbit that contains both real and complex embeddings, and let $x\in C\cap O$. Then $x$ is adjacent
to every embedding of $O$ outside its colour class in $\Gamma[C\cap O]$. Hence $\Pi(x)$ is this colour class.
\end{lemma}

\begin{proof}
We may assume that $x\in P_O$.

\smallskip\noindent\emph{Step 1: every real embedding of $O$ is adjacent to every complex embedding of $O$.} Let
$r\in R\cap O$. By Lemma~\ref{lem:connected}, $r$ has a complex neighbour, and $\Gamma[C\cap O]$ is connected. So it
suffices to show: if $r\sim a$ and $a\sim b$ with $a,b\in C$, then $r\sim b$. Let $d_1$ interchange $a$ and $r$, and
$d$ interchange $a$ and $b$. By Lemma~\ref{lem:fix}, $d$ fixes $r$. Hence $dd_1d$ interchanges $da=b$ and $dr=r$.

\smallskip\noindent\emph{Step 2: $x$ is adjacent to every $w\in Q_O$.} If $w=cx$ there is nothing to prove, so let
$w\ne cx$. Choose $r\in R\cap O$. By Step~1 there are $d_x,d_w\in\cD$ with $d_x\colon x\leftrightarrow r$ and
$d_w\colon w\leftrightarrow r$. If $d_w$ fixes $x$, then $d_wd_xd_w$ interchanges $x$ and $w$, and we are done.
Suppose that $d_w$ moves $x$. The embedding $d_wx$ cannot be complex: otherwise Lemma~\ref{lem:fix} would say that
$d_w$ fixes $r$, which is adjacent to $x$; but $d_wr=w$. Thus $d_w$ interchanges the complex embedding $x$ with a
real embedding, and by Lemma~\ref{lem:triangle} it fixes $cx$. The conjugate $e=cd_xc$ interchanges $r$ and $cx$,
so $d_wed_w$ interchanges $d_wr=w$ and $d_w(cx)=cx$. But $w$ and $cx$ both lie in $Q_O$, and two embeddings of the
same colour are not adjacent. So $d_w$ fixes $x$ after all.

\smallskip
By Steps~1 and~2, $\Pi(x)\subseteq P_O$. Conversely, $P_O\subseteq\Pi(x)$, since two embeddings of the same colour
are not adjacent.
\end{proof}

\section{The canonical subfield and the proof of the list}\label{sec:canonical}

One more notion from graph theory is needed. A graph is \emph{complete multipartite} if its vertex set is
partitioned into nonempty sets, the \emph{parts}, such that two vertices are adjacent exactly when they lie in
different parts. With two parts one speaks of a \emph{complete bipartite} graph. For instance,
Lemma~\ref{lem:parts} says that $\Gamma[C\cap O]$ is complete bipartite with parts $P_O$ and $Q_O$, for
every $N$-orbit $O$ that contains real and complex embeddings; no edge joins two embeddings of the same colour
class. If $O$ contains no real embeddings, then $\Gamma[C\cap O]=\Gamma[O]$ need not be complete bipartite: in
Example~\ref{ex:cube} it is the graph of a cube.

Assume $(\star)$. Let $d\in\cD$, and let $O$ be an $N$-orbit that meets $C(d)$. By Lemma~\ref{lem:connected},
applied to $d$, the graph $\Gamma[C(d)\cap O]$ is connected, and by $(\star)$ for $d$ it is bipartite; so it has
exactly one proper $2$-colouring, up to exchanging the colours. We call its two colour classes the \emph{halves}
of $d$ in $O$. Thus the halves of $c$ are the sets $P_O$ and $Q_O$. An element $\gamma\in G$ maps the halves of
$d$ in $O$ to the halves of $\gamma d\gamma^{-1}$ in $\gamma O$, since it is an automorphism of $\Gamma$ and maps
$C(d)$ to $C(\gamma d\gamma^{-1})$.

\begin{proposition}\label{prop:partition}
Assume $(\star)$. The halves of all $d\in\cD$, in all $N$-orbits, form a $G$-invariant partition $\fB$ of $\Sigma$
which refines the partition into $N$-orbits. The blocks of $\fB$ contained in $C$ are the sets $P_O$ and $Q_O$,
for the $N$-orbits $O$ meeting $C$, and all other blocks consist of real embeddings. Moreover, for every $N$-orbit
$O$, either $\Gamma[O]$ is bipartite and the blocks in $O$ are its two colour classes, or $\Gamma[O]$ is the
complete multipartite graph whose parts are the blocks in $O$.
\end{proposition}

\begin{proof}
By Lemma~\ref{lem:components}(a), every embedding is moved by some $d\in\cD$, and hence lies in a half of $d$. Let
$X$ and $Y$ be halves of $d$ and $d'$ in the same $N$-orbit $O$, and let $\rho\in X\cap Y$. We show that $X=Y$.

Suppose first that $\Gamma[O]$ is bipartite. It is connected by Lemma~\ref{lem:components}(b), and $d$ is an
automorphism of $\Gamma[O]$ which maps $\rho$ to a neighbour. So $d$ exchanges the two colour classes of
$\Gamma[O]$ and moves every embedding of $O$. Hence $C(d)\cap O=O$, and the halves of $d$ in $O$ are the colour
classes of $\Gamma[O]$. The same holds for $d'$, so $X=Y$.

Suppose now that $\Gamma[O]$ is not bipartite. Then $d$ fixes some embedding of $O$: otherwise $O\subseteq C(d)$,
and $\Gamma[O]=\Gamma[C(d)\cap O]$ would be bipartite by $(\star)$ for $d$. So $O$ contains embeddings fixed by $d$
and embeddings moved by $d$, and Lemma~\ref{lem:parts}, applied to $d$, shows that $X=\Pi(\rho)$. In the same way
$Y=\Pi(\rho)$, so $X=Y$. As every embedding of $O$ lies in a half, the blocks in $O$ are the sets $\Pi(\rho)$; by
the definition of $\Pi$, two embeddings of $O$ are adjacent exactly when they lie in different blocks.

Thus the halves form a partition $\fB$ of $\Sigma$. It refines the partition into $N$-orbits, and it is
$G$-invariant because $G$ maps halves to halves. If $O$ meets $C$, the halves $P_O$ and $Q_O$ of $c$ cover
$C\cap O$, so the other blocks in $O$ consist of real embeddings. If $O$ does not meet $C$, then $O\subseteq R$.
\end{proof}

\begin{definition}\label{def:Kc}
Assume $(\star)$. The \emph{canonical subfield} $K_c$ is the subfield of $K$ whose blocks form the partition
$\fB$ of Proposition~\ref{prop:partition}.
\end{definition}

It exists because $\fB$ is $G$-invariant (Section~\ref{subsec:subfields}). Since $\fB$ refines the
$N$-orbits, which are the blocks over $K^{+}$, we have $K^{+}\subseteq K_c$. By construction $K_c$ depends only
on $G$ and $c$, not on $S$.

\begin{proof}[Proof of Theorem~\ref{thm:list}, (a)$\Rightarrow$(b)]
Let $T_{K,S}$ be quasi-abelian. By Proposition~\ref{prop:colouring}, $(\star)$ holds, and $S\cap C$ is a union of
colour classes of $\Gamma[C]$; by Lemma~\ref{lem:connected}, $S\cap C\cap O$ is $P_O$ or $Q_O$ for each $N$-orbit
$O$ meeting $C$. We show that $S$ is lifted from $K_c$ and that $(K_c,K^{+})$ is admissible.

Every block of $\fB$ lies in $S$ or in $\Sigma\setminus S$: the real blocks lie in $R\subseteq S$, and the
others are the sets $P_O,Q_O$. A block mapped to itself by $c$ cannot be $P_O$ or $Q_O$, so it consists of real
embeddings. Hence $S$ is lifted from $K_c$. The embeddings of $K_c$ above a place of $K^{+}$ are the blocks inside
the corresponding $N$-orbit $O$, and the complex ones among them are those moved by $c$: the pair $P_O,Q_O$ if $O$
meets $C$, none otherwise. Thus there is at most one complex place of $K_c$ above each place of $K^{+}$, that is,
$(K_c,K^{+})$ is admissible.
\end{proof}

\begin{proof}[Proof of the last part of Theorem~\ref{thm:list}]
Suppose that some $T_{K,S}$ is quasi-abelian. Then $(\star)$ holds by Proposition~\ref{prop:colouring}, and we have
just seen that (b) holds with $K'=K_c$ for every quasi-abelian $T_{K,S}$. Now let $K'$ be a subfield for which (b)
holds, for some $S$; we show that $K'=K_c$. Since $K^{+}\subseteq K'$, the blocks over $K'$ refine the $N$-orbits.
A block over $K'$ that $c$ maps to itself corresponds to a real place of $K'$, so by (L1) it consists of real
embeddings, and $c$ fixes it pointwise. Since the partition into blocks over $K'$ is $G$-invariant, every
$d\in\cD$ likewise fixes pointwise each block over $K'$ that it maps to itself; in particular, no $d\in\cD$
interchanges two embeddings of the same block. Now let $O$ be an $N$-orbit meeting $C$. By admissibility, $c$
exchanges at most one pair of blocks inside $O$, so $C\cap O=B\sqcup cB$ for a block $B\ne cB$. These two blocks
form a proper $2$-colouring of $\Gamma[C\cap O]$, which is connected by Lemma~\ref{lem:connected}; hence
$\{B,cB\}=\{P_O,Q_O\}$. Thus the halves of $c$ are blocks over $K'$, and by $G$-invariance so are the halves of
every $d\in\cD$. Since the halves form a partition of $\Sigma$ (Proposition~\ref{prop:partition}), the blocks over
$K'$ are exactly those of $\fB$, and $K'=K_c$.

In particular, the quasi-abelian $T_{K,S}$ are exactly those with $S$ lifted from $K_c$: every quasi-abelian $S$ is
lifted from $K_c$, and conversely every $S$ lifted from $K_c$ is quasi-abelian, by the implication
(c)$\Rightarrow$(a) applied to the admissible pair $(K_c,K^{+})$. The systems lifted from $K_c$ are in bijection
with the systems of representatives of $K_c$, and there are $2^{r_2(K_c)}$ of those. This agrees with
Proposition~\ref{prop:colouring}: $r_2(K_c)$ is the number of $N$-orbits meeting $C$, which by
Lemma~\ref{lem:connected} is the number of components of $\Gamma[C]$.
\end{proof}

Together with Dodson's theorem, the list describes every quasi-abelian $T_{K,S}$ up to isogeny.

\begin{corollary}\label{cor:isogeny}
Let $T_{K,S}$ be quasi-abelian, put $m=[K:K_c]$ and $S_c=\{\iota\in\Sigma(K_c):\ \Sigma_\iota\subseteq S\}$. Then
$T_{K,S}$ is isogenous to $(T_{K_c,S_c})^{m}$, and isomorphic to it if $\OK$ is a free $\mathcal O_{K_c}$-module.
The factor $T_{K_c,S_c}$ is quasi-abelian by the trace-form criterion, applied to $K_c$ and its totally real
subfield $K^{+}$; in fact $T_{K_c,S'}$ is quasi-abelian for every system of representatives $S'$ of $K_c$.
\end{corollary}

\begin{proof}
By Theorem~\ref{thm:list}, $S$ is lifted from $K_c$, and $S_c$ is the system of representatives of $K_c$ from
which it is lifted. Theorem~\ref{thm:lifting}(a) gives the isogeny and, in the free case, the isomorphism.
Since $(K_c,K^{+})$ is admissible, Theorem~\ref{thm:trace} applies to the field $K_c$ and its subfield $K^{+}$.
\end{proof}

In other words, up to isogeny every quasi-abelian $T_{K,S}$ is a power of a group to which the trace-form
criterion applies directly, with $K^{+}$ as the totally real subfield. When $K_c=K$, as in Examples~\ref{ex:deg9}
and~\ref{ex:triangles}, the corollary says nothing new. In Example~\ref{ex:octahedron}, $K_c$ is the cubic field
$L$ with $r_2(L)=1$, and the four-dimensional group $T_{K,S_+}$ is isogenous to the square of the
two-dimensional quasi-abelian group $T_{L,S'}$, where $S'$ is the system of $L$ from which $S_+$ is lifted.

\section{Consequences}\label{sec:consequences}

\subsection{All systems of representatives at once}
The trace-form criterion gives quasi-abelian groups for all $S$ simultaneously. It is the only way this can
happen.

\begin{corollary}\label{cor:uniform}
The following are equivalent.
\begin{enumerate}[(a),leftmargin=0.9cm,itemsep=2pt]
\item $T_{K,S}$ is quasi-abelian for every system of representatives $S$.
\item $(K,K^{+})$ is admissible: above each place of $K^{+}$ there is at most one complex place of $K$.
\item $(K,K_0)$ is admissible for some totally real subfield $K_0$.
\end{enumerate}
\end{corollary}

\begin{proof}
(b)$\Rightarrow$(c) is clear, and (c)$\Rightarrow$(a) is Theorem~\ref{thm:trace}. Assume (a). By
Theorem~\ref{thm:list} all $2^{r_2}$ systems are lifted from $K_c$, and there are $2^{r_2(K_c)}$ of those; so
$r_2(K_c)=r_2$. Let $O$ be an $N$-orbit meeting $C$. Above the corresponding place of $K^{+}$ there are
$|P_O|$ complex places of $K$ and one of $K_c$. Summing over these orbits gives $r_2=\sum_O|P_O|$ and
$r_2(K_c)=\sum_O1$. Hence $|P_O|=1$ for every $O$, which is (b).
\end{proof}

For a single $S$ the trace-form criterion is not enough, and lifting is needed: see
Example~\ref{ex:octahedron}.

\subsection{The compact case}
For $r_1=0$, $T_{K,S}$ is a complex torus, which we denote $A_S$. Here $F=\C^{g}$, an ample Riemann form is a
Riemann form, and by Riemann's theorem $A_S$ is an abelian variety if and only if it is quasi-abelian in our
sense. Recall that $K'$ is a \emph{CM field} if it is a totally imaginary quadratic extension of a totally real
field, and that a system of representatives of a CM field is called a \emph{CM type}.

\begin{theorem}[Shimura--Taniyama~{\cite{ST}}; see also~{\cite[Ch.~1]{La}}]\label{thm:ST}
Let $K$ be totally imaginary. Then $A_S$ is an abelian variety if and only if $S$ is lifted from a CM type of
a CM subfield $K'\subseteq K$ (possibly $K'=K$). In that case $A_S$ is isogenous to $(A_{S'})^{[K:K']}$.
\end{theorem}

\begin{proof}
Let $S$ be lifted from a CM type $S'$ of a CM field $K'$, and let $K'_0$ be the totally real subfield with
$[K':K'_0]=2$. Above each place of $K'_0$ lies exactly one place of $K'$, and it is complex. Hence $(K',K'_0)$
is admissible, and $A_S$ is an abelian variety by Theorem~\ref{thm:list}. The isogeny is
Theorem~\ref{thm:lifting}(a).

Conversely, let $A_S$ be an abelian variety. By Theorem~\ref{thm:list}, $S$ is lifted from $K_c$, and
$(K_c,K^{+})$ is admissible. Since $r_1=0$, every $N$-orbit $O$ lies in $C$, so by Proposition~\ref{prop:partition}
the blocks of $K_c$ inside $O$ are $P_O$ and $Q_O$. Thus above each embedding of $K^{+}$ there are exactly two embeddings of $K_c$, both complex. So
$[K_c:K^{+}]=2$ and $K_c$ is totally imaginary: $K_c$ is a CM field, $K^{+}$ is its totally real subfield, and
$S'=\{\iota:\Sigma_\iota\subseteq S\}$ is a CM type of $K_c$.
\end{proof}

The converse can also be obtained from the orbit criterion directly, by counting.

\begin{proof}[Second proof of the converse]
Let $A_S$ be an abelian variety, so that condition (c) of Theorem~\ref{thm:main} holds. Since $r_1=0$, $c$ has
no fixed point on $\Sigma$, and neither has any conjugate $c'$. So (c) says $c'(S)\cap S=\emptyset$. As
$|S|=|\Sigma|/2$, this implies
\[
  c'(S)=\Sigma\setminus S=c(S)\qquad\text{for every complex conjugation }c' .
\]
This counting step is the heart of the matter. It shows that every product $c\,c'$ preserves $S$. Let $N_0$ be
the subgroup generated by these products. It is normal, since
$\gamma(c\,c')\gamma^{-1}=(c\cdot\gamma c\gamma^{-1})^{-1}(c\cdot\gamma c'\gamma^{-1})$. It preserves $S$, so it does
not contain $c$, and $N=N_0\cup cN_0$ because $c'=c\,(c\,c')$.

The $N_0$-orbits are the blocks over a subfield $K'$. As $N_0$ preserves $S$, the set $S$ is a union of
$N_0$-orbits. Since $c$ normalizes $N_0$, it maps each $N_0$-orbit $\Omega$ to the $N_0$-orbit $c\Omega$, which
lies on the other side of $S$; so $c$ fixes no block, and $K'$ is totally imaginary. The $N$-orbits are the sets
$\Omega\cup c\Omega$. By Lemma~\ref{lem:N} they are the blocks over $K^{+}$; hence $K^{+}\subseteq K'$ and
$[K':K^{+}]=2$, and $K'$ is a CM field. Finally $S$ is lifted from $S'=\{\iota:\Sigma_\iota\subseteq S\}$, which
contains exactly one of $\iota$ and $c\iota$ for every $\iota$: a CM type of $K'$.
\end{proof}

In both proofs the CM field is $K_c$, which depends only on $G$ and $c$. (In the second proof, every
$N$-orbit $O$ consists of two $N_0$-orbits, and $S\cap O$ is one of them; it is also one of $P_O$, $Q_O$.) No
endomorphism algebra and no Rosati involution enter the argument; they are needed for the finer statements of Shimura--Taniyama, such as
the description of $\operatorname{End}^0(A_S)$, which the criterion does not address. The counting step is
exactly what is lost when $r_1\ge1$: then $|S|<|\Sigma|/2$, a complex conjugation may fix embeddings in $S$ or
move them to real embeddings, and the products $c\,c'$ need not preserve $S$. This is why the general case
needs the finer analysis of Sections~\ref{sec:graph}--\ref{sec:canonical}, where the CM field is replaced by
the block system $\fB$.

\subsection{Fields without intermediate subfields}

\begin{corollary}\label{cor:primitive}
Suppose that $K$ has no subfields other than $\Q$ and $K$, or equivalently that $G$ acts primitively on
$\Sigma$. Then $T_{K,S}$ is quasi-abelian if and only if $r_2=1$.
This applies in particular when $[K:\Q]$ is prime.
\end{corollary}

\begin{proof}
For $r_2=1$ this is Corollary~\ref{cor:easy}(i). Let $T_{K,S}$ be quasi-abelian. By Theorem~\ref{thm:list},
$S$ is lifted from a subfield $K'$ with $(K',K^{+})$ admissible, and $K'$ is $\Q$ or $K$. It is not $\Q$: the
only place of $\Q$ is real, and the complex places of $K$ lie above it, which (L1) excludes. So $K'=K$, and
$(K,K^{+})$ is admissible. As $r_2\ge1$, $K$ is not totally real, so $K^{+}=\Q$, and admissibility says that $K$
has at most one complex place.
\end{proof}

In prime degree there is also a direct argument from the orbit criterion, which shows in addition that every
orbit of pairs contains a diagonal pair.

\begin{remark}\label{rem:prime}
Let $[K:\Q]=p$ be an odd prime and $r_2\ge2$. A transitive group of prime degree $p$ is $2$-transitive or
solvable, and a solvable one is conjugate to a subgroup of the affine group $\mathrm{AGL}(1,p)$
(Burnside and Galois; see for instance~\cite{DM}). In the first case Corollary~\ref{cor:easy}(ii) applies. In
the second, identify $\Sigma$ with $\mathbb F_p$ so that $G=\{x\mapsto ax+b:\ a\in U,\ b\in\mathbb F_p\}$ for a
subgroup $U\le\mathbb F_p^{\times}$; $G$ contains all translations because it contains an element of order $p$.
An element $x\mapsto ax+b$ of order $2$ has $a^{2}=1$, and $x\mapsto x+b$ has order $p$ unless $b=0$; so $c$
has the form $x\mapsto -x+b$. Translating
coordinates we may assume $c(x)=-x$; then $-1\in U$, and the only real embedding is $0$. A pair $\{x,y\}$ is
moved by translations to $\{0,y-x\}$ and then by $U$, so its $G$-orbit is determined by the coset $(y-x)U$ in
$\mathbb F_p^{\times}/U$; since $-1\in U$, this does not depend on the order of $x$ and $y$. The diagonal pair
$\{k,-k\}$ lies in the orbit of the coset $2kU$, and every coset $aU$ has this form, with $k=a/2$. Hence every
orbit contains a diagonal pair. Since $r_2\ge2$ there is a forbidden pair, and it shares its orbit with a
diagonal pair. Condition (b) of Theorem~\ref{thm:main} does not hold.
\end{remark}

\section{Meromorphic functions}\label{sec:meromorphic}

The criterion decides more than it seems to. When its conditions are not satisfied, the only meromorphic
functions on $T_{K,S}$ are the constants.

\begin{theorem}\label{thm:mero}
Let $r_2\ge1$ and $T=T_{K,S}$.
\begin{enumerate}[(a),leftmargin=0.9cm,itemsep=2pt]
\item If $r_1\ge1$, every holomorphic function on $T$ is constant: $T$ is a Cousin group.
\item Either $T$ is quasi-abelian, or every meromorphic function on $T$ is constant.
\end{enumerate}
For $r_1=0$, part (b) says that the complex torus $T$ is an abelian variety or has algebraic dimension~$0$.
In particular, every meromorphic function on $T$ is constant if and only if the equivalent conditions of
Theorems~\ref{thm:main} and~\ref{thm:list} are not satisfied.
\end{theorem}

Part (a) is known. It goes back to Andreotti and Gherardelli~\cite[\S4, Prop.~6]{AG73}, and proofs can be found
in~\cite[Lemma~2.4]{OT} and~\cite[Prop.~4.1]{AbeFM}. We include a short proof, because the same argument
gives (b).

We use two classical reductions. For a connected complex abelian Lie group $X$ there is a largest connected
closed complex subgroup $H_0$ on whose cosets every holomorphic function is constant, and
$X/H_0\cong\C^{a}\times(\C^{\times})^{b}$ (Morimoto~\cite{Mor}; see also~\cite{AK}). Suppose now that every
holomorphic function on $X$ is constant. Then there is a connected closed complex subgroup $Z_0$ of $X$ such
that $X/Z_0$ is quasi-abelian or trivial, and every meromorphic function on $X$ is the pull-back of a
meromorphic function on $X/Z_0$. This is the meromorphic reduction, or quasi-abelian reduction, of Capocasa and
Catanese~\cite[Def.~2.18 and Thm.~2.20]{CC}; see also~\cite{AK}. For a complex torus it is the classical
algebraic reduction.

Both subgroups are maximal for their property. If every meromorphic function on $X$ is constant on the cosets
of a connected subgroup $W$, then the image of $W$ in $X/Z_0$ is trivial, because a quasi-abelian group is
meromorphically separable~\cite[Main Thm.~3.2.21]{AK}; so $W\subseteq Z_0$. In the same way, holomorphic
functions separate the points of $\C^{a}\times(\C^{\times})^{b}$. Now let $\varphi$ be a surjective
holomorphic endomorphism of $X$ with finite kernel. If $f$ is meromorphic on $X$, so is $f\circ\varphi$, which is
therefore constant on the cosets of $Z_0$. Hence $f$ is constant on the cosets of $\varphi(Z_0)$, and
$\varphi(Z_0)\subseteq Z_0$. Both groups are connected, closed and of the same dimension, so
$\varphi(Z_0)=Z_0$. The same argument gives $\varphi(H_0)=H_0$.

Choose a primitive element of $K$ and multiply it by a positive integer so that it lies in $\OK$. Call the result
$\alpha$; thus $\alpha\in\OK$ and $K=\Q(\alpha)$. The diagonal matrix $A_\alpha=\mathrm{diag}(\rho(\alpha))_{\rho\in S}$
satisfies $A_\alpha u(a)=u(\alpha a)$, so $A_\alpha\Lambda\subseteq\Lambda$. Since $\alpha\ne0$, it induces a
surjective endomorphism $\varphi_\alpha$ of $T$ with finite kernel. Its eigenvalues $\rho(\alpha)$, $\rho\in S$, are
pairwise distinct, because distinct embeddings of $K=\Q(\alpha)$ differ on $\alpha$.

\begin{lemma}[rigidity]\label{lem:rigidity}
Let $W\subseteq T$ be a connected closed complex subgroup with $\varphi_\alpha(W)=W$. Then $W=0$ or $W=T$.
\end{lemma}

\begin{proof}
Let $V\subseteq\C^{g}$ be the Lie algebra of $W$, so that $W=\pi(V)$ for the projection $\pi\colon\C^{g}\to T$,
and $A_\alpha V=V$. Write $\epsilon_\rho$, $\rho\in S$, for the unit vectors of $\C^{g}$. Since the
eigenvalues of $A_\alpha$ are distinct, $V$ is spanned by some of the $\epsilon_\rho$.

The group $\Lambda\cap V$ is mapped into itself by $A_\alpha$, so it corresponds to a $\Z[\alpha]$-submodule of
$\OK$. If it contains $u(a)$ with $a\ne0$, it contains $u(a\alpha^{k})$ for all $k\ge0$, and
$a,a\alpha,\dots,a\alpha^{\delta-1}$ are linearly independent over $\Q$. Hence $\Lambda\cap V$ has rank
$0$ or $\delta$.

If the rank is $\delta$, then $\Lambda\cap V$ has finite index in $\Lambda$, so $V$ contains $\R\Lambda$, which
spans $\C^{g}$ over $\C$; hence $V=\C^{g}$ and $W=T$.

Suppose the rank is $0$ but $V\ne0$, say $\epsilon_\rho\in V$. Then $\epsilon_\rho\in\R\Lambda$, so
$V_0=V\cap\R\Lambda$ is a
nonzero real subspace. Since $W$ is closed, $V+\Lambda=\pi^{-1}(W)$ is closed, hence so is
$(V+\Lambda)\cap\R\Lambda=V_0+\Lambda$, and so is its image in $\R\Lambda/V_0$. That image is a countable closed
subgroup of a real vector space of dimension $<\delta$. A closed subgroup that is not discrete contains a line,
which is uncountable. So the image is discrete, and its rank is $<\delta$. But it is isomorphic to
$\Lambda$, since $\Lambda\cap V_0=0$, and $\Lambda$ has rank $\delta$. This contradiction shows $V=0$.
\end{proof}

\begin{proof}[Proof of Theorem~\ref{thm:mero}]
(a) Let $r_1\ge1$, and let $H_0$ be the subgroup of the holomorphic reduction. It is preserved by $\varphi_\alpha$, so
$H_0=0$ or $H_0=T$ by Lemma~\ref{lem:rigidity}. If $H_0=0$, then $T\cong\C^{a}\times(\C^{\times})^{b}$ with
$a+b=g$, and its fundamental group $\Z^{b}$ has rank at most $g$. But the fundamental group of $T$ is $\Lambda$,
of rank $\delta=g+r_2>g$. Hence $H_0=T$.

(b) Every holomorphic function on $T$ is constant: by (a) if $r_1\ge1$, and because $T$ is compact if
$r_1=0$. So $T$ has a meromorphic reduction $Z_0$, again preserved by $\varphi_\alpha$. If $Z_0=0$, then
$T=T/Z_0$ is quasi-abelian. If $Z_0=T$, every meromorphic function is the pull-back of a function on a point,
and so it is constant. The last assertion follows from Theorems~\ref{thm:main} and~\ref{thm:list}, since a
quasi-abelian group of positive dimension has nonconstant meromorphic functions~\cite{CC,Abe89}.
\end{proof}

For the groups $T_{K,S}$ there is thus no intermediate case between being quasi-abelian and having only
constant meromorphic functions. When $r_1\ge1$, one may even take for $\alpha$ a unit, by a theorem of
Za\"imi, Bertin and Aljouiee~\cite[Thm.~1.4]{ZBA}; then $\varphi_\alpha$ is an automorphism. This is not needed
here, and it is not available in general when $r_1=0$. Kazama and Umeno~\cite{KU} proved that a toroidal group whose N\'eron--Severi
group vanishes has only constant meromorphic functions, and Umeno used this for $\Q(2^{1/5})$~\cite{Um17}.
Theorem~\ref{thm:mero} needs no such hypothesis: for $S_-$ in Example~\ref{ex:octahedron} and in the dihedral
case of Example~\ref{ex:quintic}, the N\'eron--Severi group has positive rank, and still every meromorphic
function is constant. A general Cousin group can have a proper meromorphic reduction, for instance
a product of a quasi-abelian group with one that has only constant meromorphic functions.

\section{Examples and remarks}\label{sec:examples}

In all examples the decision comes from the block structure over a subfield or from the shape of the orbits.
We describe the conjugation graph $\Gamma$ where it is small.

\begin{example}[degree $9$, signature $(5,2)$]\label{ex:deg9}
Let $\mu^{3}-4\mu-1=0$ and $K=\Q(\alpha)$ with $\alpha^{3}-\alpha=\mu$. The cubic $\Q(\mu)$ is totally real; over
its three real places the polynomial $x^{3}-x-\mu$ has $3$, $1$ and $1$ real roots, so $b=(0,1,1)$. By
Theorem~\ref{thm:trace}, $T_{K,S}$ is quasi-abelian for all four systems $S$. Here $K^{+}=\Q(\mu)$, since $K$ is not
totally real and no field lies strictly between $\Q(\mu)$ and $K$; and $K_c=K$. The same holds for $K=\Q(\eta)$ with
$\eta^{3}+4\beta\eta+4=0$ over $\Q(\beta)=\Q(\zeta_9)^{+}$, $\beta^{3}-3\beta+1=0$, where $b=(1,1,0)$.
\end{example}

\begin{example}[two triangles]\label{ex:triangles}
Let $K=\Q(\sqrt2,\sqrt[3]3)$, of signature $(2,2)$. Over each real place of $K^{+}=\Q(\sqrt2)$ lie one real and one
complex place, so all four $S$ give quasi-abelian groups by Theorem~\ref{thm:trace}. Two of them are lifted from
$\Q(\sqrt[3]3)$, which has $r_2=1$; the other two are not lifted from any proper subfield. The complex
conjugations act on the three cube roots of $3$ as the three transpositions and fix $\sqrt2$. Hence $\Gamma$
consists of two triangles, one for each embedding of $\Q(\sqrt2)$, each with one real and two complex embeddings.
The parts $\Pi(\rho)$ consist of single embeddings, and $K_c=K$. This is an instance of the situation of Andreotti and
Gherardelli: a cubic extension of a totally real field. By Proposition~\ref{prop:NS} the N\'eron--Severi group
has rank $12$ for the two systems lifted from $\Q(\sqrt[3]3)$ and rank $9$ for the other two.
\end{example}

\begin{figure}[ht]
\centering
\begin{tikzpicture}[scale=1.05,
  v/.style={circle,draw,thick,fill=white,minimum size=7mm,inner sep=0pt,font=\small},
  e/.style={thick,black!70}]
\node[v] (s1)  at (90:2.3)   {$\sigma_1$};
\node[v] (t)   at (210:2.3)  {$\theta$};
\node[v] (tb)  at (330:2.3)  {$\bar\theta$};
\node[v] (s2)  at (270:0.62) {$\sigma_2$};
\node[v] (mt)  at (30:0.62)  {$-\theta$};
\node[v] (mtb) at (150:0.62) {$-\bar\theta$};
\foreach \a/\b in {s1/t,t/tb,tb/s1,s2/mt,mt/mtb,mtb/s2,s1/mt,s1/mtb,t/s2,t/mtb,tb/s2,tb/mt}
  \draw[e] (\a) -- (\b);
\end{tikzpicture}
\caption{The conjugation graph $\Gamma=K_{2,2,2}$ of Example~\ref{ex:octahedron}. Opposite vertices form the
three blocks over $L$; they are the only pairs not joined by an edge. The complex embeddings span the
$4$-cycle $\theta,\bar\theta,-\theta,-\bar\theta$.}
\label{fig:octahedron}
\end{figure}
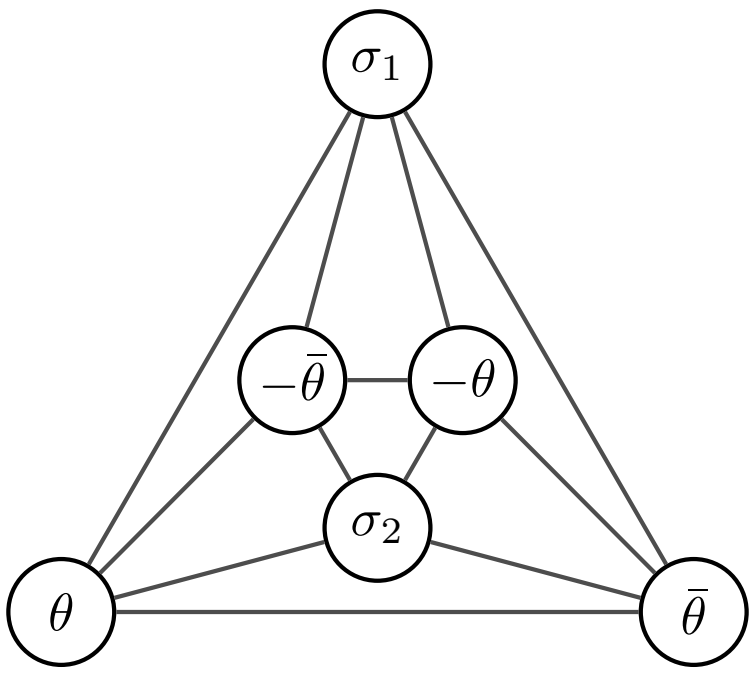

\begin{example}[the octahedron]\label{ex:octahedron}
Let $K=\Q(\theta)$ with $\theta^{6}-6\theta^{4}+12\theta^{2}-10=0$, of signature $(2,2)$. Here
$\theta^{2}=2+\lambda$ with $\lambda^{3}=2$, and $G=C_2\wr S_3$, of order $48$. Fix a non-real root of the
polynomial and call it $\theta$ again. Besides the real embeddings $\sigma_1$ and $\sigma_2$ there are four
complex ones, which we name after the image of the generator: $\theta$, $-\theta$, $\bar\theta$ and
$-\bar\theta$. Over $L=\Q(\lambda)$, which has one real and one complex place, the blocks are
$\{\sigma_1,\sigma_2\}$, $\{\pm\theta\}$ and $\{\pm\bar\theta\}$, and $c$ interchanges the last two. The only proper subfields are $\Q$ and $L$, so $K^{+}=\Q$ and there is one $N$-orbit.
The graph $\Gamma$ is the graph of the octahedron, with its six vertices and twelve edges: the complete
multipartite graph $K_{2,2,2}$ whose parts are the three blocks over $L$. Each embedding is adjacent to the four
embeddings outside its block, and the blocks are the three pairs of opposite vertices; see Figure~\ref{fig:octahedron}. Hence $K_c=L$, and exactly two systems give quasi-abelian groups, the two systems lifted from $L$.
For $S_+=\{\sigma_1,\sigma_2,\theta,-\theta\}$ the group is quasi-abelian. For $S_-=\{\sigma_1,\sigma_2,\theta,-\bar\theta\}$
it is not: the element $\gamma\in G$ that changes the sign on $\{\pm\bar\theta\}$ and fixes the other embeddings
gives a complex conjugation $\gamma c\gamma^{-1}$ that interchanges $\theta$ and $-\bar\theta$. By
Theorem~\ref{thm:mero}, $T_{K,S_-}$ has only constant meromorphic functions, although its N\'eron--Severi group
has rank $3$ by Proposition~\ref{prop:NS}: the three pairs inside the blocks over $L$ form an orbit without
forbidden pairs. So the answer depends on $S$, not only
on $K$. Since $K^{+}=\Q$ and $r_2=2$, no pair $(K,K_0)$ is admissible: the quasi-abelian group $T_{K,S_+}$ comes
from lifting, not from the trace-form criterion on $K$.
\end{example}

\begin{example}[the cube]\label{ex:cube}
Let $\beta_1,\beta_2,\beta_3$ be the roots of $t^{3}-3t+1$, so that $K_0=\Q(\beta_1)=\Q(\zeta_9)^{+}$ is a cyclic
cubic field. Exactly one $\beta_j$ is negative, and $\beta_1\beta_2\beta_3=-1$. Choose square roots $r_j$ of
$\beta_j$ and put $\theta=r_1+r_2+r_3$ and $u=r_1r_2r_3$, so that $u^{2}=-1$. Since $\sum_j\beta_j=0$ and
$\sum_{j<k}\beta_j\beta_k=-3$, one finds $\theta^{4}=-12+8u\theta$, and hence
\[
  \theta^{8}+24\theta^{4}+64\theta^{2}+144=0 .
\]
Let $K=\Q(\theta)$. No product of a nonempty subset of $\{\beta_1,\beta_2,\beta_3\}$ is a square in $K_0$, since
each such product is negative at some real place of $K_0$. So $\Kt=K_0(r_1,r_2,r_3)$ has degree $24$, and
$G=(\Z/2)^{3}\rtimes C_3$: the normal subgroup $(\Z/2)^{3}=\Gal(\Kt/K_0)$ changes the signs of the $r_j$
independently, and $C_3$ permutes them cyclically. The eight numbers $\pm r_1\pm r_2\pm r_3$ are distinct, so
$[K:\Q]=8$, and $K$ is the field fixed by $C_3$. It is totally imaginary and contains $u=(\theta^{4}+12)/(8\theta)$,
a square root of $-1$.

Complex conjugation $c$ changes the sign of the $r_j$ with $\beta_j<0$ and fixes the other two. Thus $\cD$
consists of the three single sign changes, $N=(\Z/2)^{3}$, $K^{+}=\Q$, and there is one $N$-orbit. As $N$ acts
simply transitively on $\Sigma$, we may identify $\Sigma$ with $(\Z/2)^{3}$, and then $\Gamma$ is the graph of a
cube. It is bipartite, with the two sets of four vertices of even and of odd weight as colour classes, but it is
not complete bipartite: two opposite vertices have different colours and are not adjacent. The colour classes
are the blocks over $\Q(u)=\Q(\iu)$, since $u$ is fixed by the sign changes of even weight and changes sign under
the others. Hence $K_c=\Q(\iu)$. Exactly two of the sixteen systems give quasi-abelian groups, which are here
abelian varieties: the two systems lifted from $\Q(\iu)$. By Corollary~\ref{cor:isogeny}, they are isogenous to
$E^{4}$, where $E=\C/\Z[\iu]$.
\end{example}

\begin{example}[quintic fields]\label{ex:quintic}
For $K=\Q(2^{1/5})$, $G$ is the Frobenius group of order $20$, which is $2$-transitive, and for the field defined
by $x^{5}-5x+12$, $G$ is dihedral of order $10$. Both have signature $(1,2)$, and in both cases $T_{K,S}$ is not
quasi-abelian for any $S$, by Corollary~\ref{cor:primitive}. In the dihedral case the embeddings are the vertices
of a pentagon; the orbits of pairs are the five sides and the five chords. One diagonal pair is a side and the
other a chord, so every forbidden pair shares its orbit with a diagonal pair, although $\NS(T_{K,S})$ has
rank $5$ by Proposition~\ref{prop:NS}. Corollary~\ref{cor:primitive} applies to every quintic field of signature
$(1,2)$. The three-dimensional groups $T_{K,S}$ of such fields are studied in~\cite[Thms.~4.2--4.4]{DF}
and~\cite[\S3.2]{DB}. None of them is quasi-abelian, whatever the choice of $S$, and by
Theorem~\ref{thm:mero} all their meromorphic functions are constant. This contains the computations
of~\cite{Um17} and~\cite[\S3.2.1]{DB}, which treat one choice of $S$ for particular quintic fields.
\end{example}

\begin{example}[$K=\Q(2^{1/9})$, signature $(1,4)$]\label{ex:deg9b}
The embeddings are $2^{1/9}\zeta_9^{k}$, $k\in\Z/9$, and $G=\Z/9\rtimes(\Z/9)^{\times}$ acts by $k\mapsto ak+b$.
There are two orbits of pairs, according to whether the difference of the two exponents is a unit or lies in
$\{3,6\}$. The diagonal pairs $\{1,8\}$ and $\{3,6\}$ lie in different orbits, so every orbit contains a diagonal
pair, and $T_{K,S}$ is not quasi-abelian for any of the $16$ systems $S$. In terms of Theorem~\ref{thm:list}: the
only proper subfield other than $\Q$ is $\Q(2^{1/3})$, and $S$ cannot be lifted from it, because above its real
place there are complex places of $K$; and $K^{+}=\Q$ while $r_2=4$.
\end{example}

\begin{remark}[what the proof uses]\label{rem:what}
Apart from the orbit criterion and Lemma~\ref{lem:N}, the proof of Theorem~\ref{thm:list} in
Sections~\ref{sec:graph}--\ref{sec:canonical} uses only that $\cD$ is a conjugacy class of involutions in a
transitive permutation group and that $\Gamma[C]$ is bipartite. It is therefore a statement about such classes:
if the elements moved by one involution of the class induce a bipartite subgraph of $\Gamma$, then each component of
$\Gamma$ is bipartite or complete multipartite, and the parts form a block system.
\end{remark}

\section*{Use of AI tools}
In accordance with current editorial policies, the author declares the following use of the AI system Claude
(Anthropic) during this work. The AI system is not an author of this article. The main idea, the project and
the essential bibliographical input come solely from the author.

\emph{Arguments.} The AI system was used to develop and to test arguments, mainly for Section~\ref{sec:NS} and for the proof of the
complete list (Sections~\ref{sec:graph}--\ref{sec:canonical}). It also proposed the present organisation of
Sections~\ref{sec:graph}--\ref{sec:canonical} and the proof of Lemma~\ref{lem:descent} in its present form. The computational claims
in the examples of Section~\ref{sec:examples} (irreducibility and signatures of the polynomials, Galois
groups, ranks of N\'eron--Severi groups) were checked with its help by computer, and the Galois group of
Example~\ref{ex:cube} was confirmed independently with Maple. The main results, their formulation and the
strategy of proof are due to the author.

\emph{Writing.} The text of Sections~\ref{sec:NS} and~\ref{sec:graph}--\ref{sec:canonical} and the examples of
Section~\ref{sec:examples} were written in part with the AI system. For the other sections the AI system was used only to revise wording and presentation. The AI system also helped
with the preparation of the \LaTeX{} source.

\emph{Citations.} The bibliographical input is essentially due to the author; the AI system was used only
occasionally to verify references. All references, including the numbering of the cited theorems, were
checked by the author against the original sources.

All statements and proofs have been checked by the author, who takes full responsibility for the content of
the article.

\end{document}